\documentclass[12pt,leqno]{amsart}
\usepackage{amssymb,amsmath,amsthm,amsfonts}
\usepackage{graphicx}
\usepackage{color}
\usepackage[dvipsnames]{xcolor}

\DeclareMathOperator{\diam}{diam}
\DeclareMathOperator{\dist}{dist}
\newcommand{\field}[1]{\mathbb{#1}}
\newcommand{\N}{\field{N}}                      
\newcommand{\R}{\field{R}}                      
\newcommand{\C}{\field{C}}                      

\newcommand{\eps}{\epsilon}

\newcommand{\loc}{{\scriptstyle{loc}}}

\newcommand{\RH}{{\scriptstyle{RH}}}

\newcommand{\Sob}{{\scriptstyle{Sob}}}
\newcommand{\reg}{{\scriptstyle{reg}}}

\newcommand{\bA}{{\mathbf A}}

\newcommand{\cI}{{\mathcal I}}

\newcommand{\bv}{{\mathbf v}}

\newcommand{\cW}{{\mathcal W}}

\newcommand{\ut}{{\underline{t}}}
\newcommand{\ot}{{\overline{t}}}

\newcommand{\Lip}{\operatorname{Lip}}

\newcommand{\bi}{{\mathbf i}}

\def\Barint_#1{\mathchoice
          {\mathop{\vrule width 6pt height 3 pt depth -2.5pt
                  \kern -8pt \intop}\nolimits_{#1}}%
          {\mathop{\vrule width 5pt height 3 pt depth -2.6pt
                  \kern -6pt \intop}\nolimits_{#1}}%
          {\mathop{\vrule width 5pt height 3 pt depth -2.6pt
                  \kern -6pt \intop}\nolimits_{#1}}%
          {\mathop{\vrule width 5pt height 3 pt depth -2.6pt
                  \kern -6pt \intop}\nolimits_{#1}}}

\theoremstyle{plain}
\newtheorem{theorem}{Theorem}
\newtheorem{corollary}[theorem]{Corollary}
\newtheorem{lemma}[theorem]{Lemma}
\newtheorem{proposition}[theorem]{Proposition}

\theoremstyle{definition}
\newtheorem{definition}[theorem]{Definition}
\newtheorem{example}[theorem]{Example}
\newtheorem{remark}[theorem]{Remark}

\newtheorem{claim}[theorem]{Claim}

\numberwithin{theorem}{section} \numberwithin{equation}{section}

\title[QC distortion of the lower dimension and applications]{Quasiconformal distortion of the lower dimension and applications}
\date{\today}
\author{Efstathios K. Chrontsios Garitsis}
\address{Department of Mathematics \\  The Ohio State University\\ 231 W 18th Ave\\ Columbus, OH 43210}
\email{chrontsios.1@osu.edu, echronts@gmail.com}
\author{Jeremy T. Tyson}
\address{Department of Mathematics \\ University of Illinois at Urbana-Champaign \\ 1409 West Green Street \\ Urbana, IL 61801}
\email{tyson@illinois.edu}


\begin{document}

\begin{abstract}
We study the distortion of the lower dimension and lower spectrum of sets under quasiconformal mappings. As an application, we verify weighted Hardy--Sobolev inequalities in Euclidean domains arising as quasiconformal images of domains with sufficiently regular boundaries. The regularity in question is phrased in terms of estimates for either the Assouad or the lower dimension of the complement. We also provide dilatation-dependent classification results for certain planar open sets.
\end{abstract}

\maketitle

\section{Introduction}

The theory of Hardy--Sobolev inequalities effectively demonstrates the close interplay between analytic inequalities and the geometry of domains. The topic is classical, originating in Hardy's investigation of Hilbert's inequality. Hardy stated the general integral inequality for $p>1$ in 1920 \cite{Hardy1920} and published a proof of its sharp form on the half-line five years later \cite{Hardy1925}. We refer to \cite{HardySurvey} for a detailed account of the history behind Hardy's inequality. The classical integral inequality is an $L^p$ bound for the averaging operator on $(0,\infty)$, which yields a weighted derivative estimate. More recent formulations of these inequalities involve control of weighted norms of functions by weighted norms of their gradients, where the weights in question involve distance to the boundary of the domain. The validity of such inequalities reflects quantitative geometric properties of the boundary of the underlying domain \cite{HardyApplLehr}, but their significance extends beyond analysis and geometry, for instance to spectral theory and mathematical physics \cite{FrankLiebSeiringer2008} or to the regularity theory of stable solutions to semilinear elliptic equations in high dimension \cite[(1.1.4)]{cabre:holder}. In recent years, the theory has been developed in metric measure spaces, where often upper gradients and Poincar\'e inequalities provide the appropriate framework \cite{HardyMetric1}, \cite{HardyMetric2}, \cite{HardyMetric3}.

Quasiconformal and quasisymmetric mapping theory provides another setting in which analytic and geometric properties are closely related. It is a field with a rich history and connections to diverse areas, such as elliptic partial differential equations \cite{aim:2qc}, nonlinear potential
theory \cite{HeinonenWeightedSobBook}, and geometric group theory \cite{mt:conf-dim-book, Bonk_Cannon_conj}. We refer to the books of V\"ais\"al\"a \cite{Vai:book} and Heinonen \cite{hei:lectures} for detailed expositions. The \textit{uniformization problem} is a central question in quasiconformal mapping theory. Roughly speaking, it asks when a space can be transformed in a controlled way into a `canonical' space with better geometric properties. The Euclidean unit interval, unit disk, and unit ball are examples of what is perceived as a canonical space. The Riemann Mapping Theorem is a signature result of this kind, but versions of this principle hold in higher Euclidean dimensions and metric spaces, the latter motivated by applications in geometric group theory and topology. A related question is to classify spaces, typically subsets of Euclidean spaces, up to quasiconformal equivalence in a quantitative way, via the dilation of the mappings in question. The celebrated paper of Gehring and V\"ais\"al\"a \cite{GeVaMonster} made a gigantic contribution to this area, and these results have been further developed by many different authors throughout the years (see \cite{BeKaiUnif}, \cite{Unif2}, \cite{Unif3}, \cite{Unif4} for a non-exhaustive list). We also refer to the surveys by Ntalampekos \cite{ntal:survey} and Bonk \cite{Bonk_QCgeom} for more information about the uniformization problem.

With an eye towards the behavior of Hardy inequalities under quasiconformal mappings (Theorem~\ref{Thm: Hardy application R2}) and the quasiconformal classification problem for  
domains with cusps (Theorem~\ref{th:ntalampekos-lower-bound}), we establish quantitative dimension distortion estimates for quasiconformal maps.
The broader program of quasiconformal dimension distortion has a long history of its own and dates back to the work of Gehring and V\"ais\"al\"a in the early 1970s \cite{GehringVais, Gehring73}. Its development is closely connected to the higher integrability of quasiconformal mappings \cite{Gehring73}. Astala's celebrated theorem \cite{Astala} provides sharp conclusions for the Sobolev regularity of planar quasiconformal maps, and thereby yields sharp bounds for the distortion of Hausdorff dimension. Kaufman's decisive contribution in \cite{Kaufman} was the introduction of an appropriate stopping-time argument in this context, which enabled extension of such distortion results to other notions of dimension such as box-counting and packing dimension. Kaufman's work inspired an avalanche of results by many different authors and the field remains very active. A partial and incomplete list of papers in this direction includes \cite{FraTyson_intermed_dim}, \cite{AC-G26}, \cite{SobSubcriticalDimDist1}, \cite{SobSubcriticalDimDist2}, \cite{OurQCspec}, \cite{SobSubcritDimDist0}, \cite{HolomSpecChron}, \cite{btw:heisenberg}, \cite{BaloghAGMS}, \cite{Chron_comp_holder_Minkowski}, and \cite{BishHakWill:QSdistortion}. A survey tracing the development of the subject from these early papers up to the most recent contributions can be found in \cite{Tys:survey}.

Different notions of dimension capture different aspects of size and local scaling. A relatively new dimension notion, at least compared to the Hausdorff dimension, is that defined by Assouad in \cite{Assouad77} (see also \cite{Assouad79, Assouad83}) which now bears his name.
More recently, Fraser and Yu introduced the Assouad spectrum \cite{fy:assouad-spectrum}, a one-parameter family of dimensions naturally interpolating between the upper box dimension and the Assouad dimension. We refer to the book by Fraser \cite{FraserBook} for a thorough introduction to the Assouad dimension and spectrum, along with applications of these notions to fractal geometry, probability, number theory, and functional analysis.
In this paper we focus on the lower (Assouad) dimension, introduced by Larman \cite{Larman1967}, and the lower spectrum, introduced by Fraser and Yu \cite{fy:assouad-spectrum}, which are essentially the dual notions to the Assouad dimension and spectrum, respectively. 
Recent work on these dimensions includes quasi-Lipschitz invariance \cite{ChenLowerSpec}, regularity properties of measures \cite{KatTros_lower_dim_spec}, and the classification of attainable lower spectra \cite{BanajiChenRutarWang2026}. One-parameter families such as the Assouad and lower spectra reside within the broader framework of {\it dimension interpolants} which has been heavily advocated by Falconer, Fraser, and their collaborators. A nice overview of the topic can be found in \cite{fraser:interpolating}.

This paper, an immediate descendant of our earlier contribution \cite{OurQCspec}, continues our program of studying the mapping properties of dimension interpolants, and identifying applications of the resulting estimates. 
Our first main result quantifies the distortion of lower dimension $\dim_L$ under quasiconformal mappings of Euclidean domains. As in \cite{OurQCspec}, the distortion bounds for the lower dimension  are governed by the {\it reverse H\"older inequality higher integrability exponent}.  
More precisely, we denote by $p_O^{\RH}(d,K)$ the supremal exponent $p$ for which the reverse H\"older inequality \eqref{eq:pnK0} holds for all cubes $Q$ contained sufficiently deep within $\Omega$ and all quasiconformal maps $f$ between domains in $\R^d$ with outer dilatation $K_O(f) \le K$. We refer to Section \ref{subsec:qc-mappings} for precise definitions and details. See also Theorem~\ref{thm:KI QC L dimension distortion} for improved but more technical bounds involving the inner dilatation.

\begin{theorem}\label{thm: QC L dimension distortion}
Let $f:\Omega\rightarrow\Omega'$ be a $K$-quasiconformal map between domains in $\R^d$, $d \ge 2$.

\medskip

\noindent (a) For any compact set $E\subset \Omega$, it holds true that
\begin{equation}\label{eq: dimL}
 \frac{(p_O^{\RH}(d,K^{d-1})-d)\,\dim_{L} E}{p_O^{\RH}(d,K^{d-1})-\dim_{L} E}\leq \dim_{L} f(E) \le \frac{p_O^{\RH}(d,K)\,\dim_{L} E}{p_O^{\RH}(d,K)\, - d + \dim_{L} E}
\end{equation}
(b) If $\Omega = \Omega' = \R^d$, then the estimate in \eqref{eq: dimL} holds true for arbitrary sets $E \subset \R^d$.
\end{theorem}

Following similar arguments to those in the proof of Theorem~\ref{thm: QC L dimension distortion}, we also establish corresponding bounds for the lower box dimension under quasiconformal mappings. In fact, we provide estimates for the distortion of lower box-counting dimension for the more general class of supercritical Sobolev mappings, complementing the result of Kaufman in \cite{Kaufman}. This in turn yields lower box dimension distortion estimates for quasiconformal mappings which involve the sharper {\it Sobolev higher integrability exponent} $p_O^{\Sob}(d,K)$. We refer to Section \ref{sec:lowerboxdim} and specifically to Theorem~\ref{thm: Sob ldimB} and Corollary \ref{cor: QC ldimB} for details.

The presentation of these dimension distortion estimates in \eqref{eq: dimL} is similar to that in Gehring's original paper \cite[Theorem 3]{Gehring73}. Many subsequent references recast these estimates in a `symmetric' framework to highlight the dependence of the estimates on the H\"older exponent of $K$-quasiconformal mappings. For lower dimension, the presentation of \eqref{eq: dimL} in this symmetric framework reads as follows:
\begin{equation}\label{eq: dim-l-1}
\left( 1 - \frac{d}{p_O^{\RH}(d,K)} \right) \left( \frac1{\dim_{L} E} - \frac1d \right)  \le \frac1{\dim_{L} f(E)} - \frac1d
\end{equation}
and
\begin{equation}\label{eq: dim-l-2}
\frac1{\dim_{L} f(E)} - \frac1d \le 
\left( 1 - \frac{d}{p_O^{\RH}(d,K^{d-1})} \right)^{-1} \left( \frac1{\dim_{L} E} - \frac1d \right).
\end{equation}
While the symmetric reformulation in \eqref{eq: dim-l-1} and \eqref{eq: dim-l-2} is merely a cosmetic change to the direct formulation in \eqref{eq: dimL}, recent contributions to quasiconformal dimension theory suggest that the `relative reciprocal dimension'
$$
\frac1{\dim E} - \frac1d, \qquad E \subset \R^d,
$$
plays a more substantial role than originally envisioned. See Remark \ref{rem:why-reciprocals} for more details.

We next present our main result on the (regularized) lower spectrum $\dim_{L,\reg}^\theta$. Recalling notation introduced in our work on quasiconformal distortion of the Assouad spectrum \cite{OurQCspec}, we set
$$
\theta(t) = \frac{1}{1+t}.
$$
For $t>0$, $K \ge 1$, and $d \ge 2$ we also set
$$
\ut = \ut(d,K) := \left(1 - \frac{d}{p_O^{\RH}(d,K)} \right)^{-1} t
$$
and
$$
\ot = \ot(d,K) := \left(1 - \frac{d}{p_O^{\RH}(d,K^{d-1})} \right) t.
$$

\begin{theorem}\label{thm: QC L spectrum distortion}
Let $f:\Omega\rightarrow\Omega'$ be a $K$-quasiconformal mapping between domains in $\R^d$, $d \ge 2$. Suppose $t\in (0,\infty)$ and $E\subset \Omega$ is compact. Then
\begin{equation}\label{eq: dim-l-spec}
\frac{(p_O^{\RH}(d,K^{d-1})-d)\,\dim_{L,\reg}^{\theta(\ot/K^{d-1})} E}
{p_O^{\RH}(d,K^{d-1})-\dim_{L,\reg}^{\theta(\ot/K^{d-1})} E}
\le \dim_{L,\reg}^{\theta(t)} f(E)
\le
\frac{p_O^{\RH}(d,K)\,\dim_{L,\reg}^{\theta(K \ut)} E}
{p_O^{\RH}(d,K)-d+\dim_{L,\reg}^{\theta(K \ut)} E}.
\end{equation}
\end{theorem}

We refer to Theorem~\ref{thm:KI QC L spectrum distortion} for once again improved bounds involving the inner dilatation.
When $d=2$ it follows from Astala's theorem \cite{Astala} that
$$
p_O^\RH(2,K) = \frac{2K}{K-1}.
$$
In this case we have
$$
\ut = Kt \qquad \mbox{and} \qquad \ot = t/K
$$
and the inequalities \eqref{eq: dimL}, respectively \eqref{eq: dim-l-spec}, read as follows:
  \begin{equation}\label{eq: dim-l-dimension-2}
\frac1K \left( \frac1{\dim_{L} E} - \frac12 \right)  \le \frac1{\dim_L f(E)} - \frac12
\le 
K \left( \frac1{\dim_{L} E} - \frac12 \right)
\end{equation}
and
\begin{equation}\label{eq: dim-spec-dimension-2}
\frac1K \left( \frac1{\dim_{L,\reg}^{\theta(K^2 t)} E} - \frac12 \right)  \le \frac1{\dim_{L,\reg}^{\theta(t)} f(E)} - \frac12
\le 
K \left( \frac1{\dim_{L,\reg}^{\theta(t/K^2)} E} - \frac12 \right).
\end{equation} 
The formulations in \eqref{eq: dim-l-dimension-2} and 
\eqref{eq: dim-spec-dimension-2} demonstrate the precise dependence of these dimension distortion bounds on the corresponding H\"older exponent of planar $K$-quasiconformal maps.

As was the case for Assouad spectrum, the value of the lower spectrum parameter $\theta$ changes under a quasiconformal mapping. However, the loss of information in the spectrum parameter is more severe for lower spectrum than it was for Assouad spectrum: in the latter case the spectrum parameter $\theta(t)$ for the image $f(E)$ was compared to $\theta(t/K)$ and $\theta(Kt)$ for the source $E$, while in the case of the lower spectrum, the spectrum parameter $\theta(t)$ for the image $f(E)$ is compared to $\theta(K\ut)$ and $\theta(K\ot)$ for the source $E$. Here $\ut$ and $\ot$ are new values which themselves depend on $t$ and also on the parameters $d$ and $K$. When $d=2$ the relevant spectrum parameters for the source $E$ in the lower spectrum case are $\theta(K^2t)$ and $\theta(t/K^2)$. This distinction arises from the rigidity of the allowed coverings employed in the lower spectrum proof. See Remark~\ref{rem: spectrum exp diff} for further details.

We now turn to the applications of Theorems~\ref{thm: QC L dimension distortion} and \ref{thm: QC L spectrum distortion}. First, we record the Hardy-Sobolev inequalities that the quasiconformal image of a Euclidean subset may satisfy. This is an application of Theorem~\ref{thm: QC L dimension distortion} and \cite[Theorem 1.2]{OurQCspec}.
\begin{theorem}\label{Thm: Hardy application R2}
Let $f:\R^2\rightarrow\R^2$ be a $K$-quasiconformal mapping, let $\Omega \subset \R^2$ be an open set, and let $q\in(1,\infty)$, $\beta<q-1$.

\medskip

\noindent (a) If $\dim_A(\R^2 \setminus \Omega) < 2 - s$,
then $f(\Omega)$ admits a $(q,\beta)$-Hardy inequality for all $q\in(1,\infty)$, $\beta<q-1$ with
$$
q-\beta
<
\frac{2s}{\,2+(K-1)(2-s)}.
$$
(b) If either $\Omega$ is bounded or $\Omega$ and $\R^2\setminus\Omega$ are both unbounded, and if $\dim_L(\R^2 \setminus \Omega) > 2 - s$, then $f(\Omega)$ admits a $(q,\beta)$-Hardy inequality for all $q\in(1,\infty)$, $\beta<q-1$ with
        $$
        q-\beta
> \frac{2Ks}{2+(K-1)s}.
        $$
\end{theorem}

We also apply our distortion results to classify various planar open sets at the level of dilatation. We highlight one such result (see Section \ref{subsec:qc-classification-planar-open-sets} for more). The `cusp domains' $\Omega_\alpha := \{ z = x + \bi y:0 <x<1, |y|<x^\alpha\}$, $\alpha>1$, are natural model examples of simply connected planar domains which are quasiconformally distinct from the unit disc. Quasiconformal equivalence of planar domains with exterior cusps and other singularities was recently considered by Ntalampekos \cite{nt:characterization}. Corollary 4.3 in \cite{nt:characterization} asserts that any two domains $\Omega_\alpha$ and $\Omega_\beta$, $1<\alpha,\beta<\infty$ are quasiconformally equivalent, 
but the dilatation constant of the mapping in question is far from optimal. As an application of Theorem \ref{thm: QC L spectrum distortion} we obtain a lower bound on the dilatation of any quasiconformal map $f:\R^2 \to \R^2$ such that $f(\Omega_\alpha) = \Omega_\beta$. 
The precise value of the lower spectrum $\dim_L^\theta(\Omega_\alpha)$ is computed in Example~\ref{ex:cusp-domains}. It is important to note that, in contrast with notions such as the Hausdorff, box-counting, or Assouad-type dimensions, the lower dimension and lower spectrum of a nonempty open subset of $\R^d$ can be strictly less than $d$. This feature, while unusual from the perspective of pure dimension theory, is essential for applications to quasiconformal equivalence as in Theorem  \ref{th:ntalampekos-lower-bound}.

 \begin{theorem}\label{th:ntalampekos-lower-bound}
 Let $f:\R^2 \to \R^2$ be a homeomorphism with $f(\Omega_\alpha) = \Omega_\beta$, $1<\alpha<\beta$. Then
 $$
 K(f)^3 + K(f)^2 \ge 2 \frac{\beta-1}{\alpha-1}.
 $$
 \end{theorem}

 An explicit lower bound can be found by solving the cubic equation. We obtain
 \begin{equation}\label{eq:Klb}
 K \ge \left( q - \frac1{27} + \sqrt{q^2 - \frac2{27}q} \right)^{1/3} +
 \left( q - \frac1{27} - \sqrt{q^2 - \frac2{27}q} \right)^{1/3} - \frac13, \qquad q = \frac{\beta-1}{\alpha-1}.
\end{equation}

We conclude the Introduction by highlighting a few of the principal differences in our current work as compared to the approach for the Assouad dimension and spectrum in \cite{OurQCspec}. If one naively tries to apply the methods in \cite{OurQCspec} that yield Assouad dimension and spectrum bounds to the lower dimension and spectrum, respectively, one quickly encounters serious obstacles. The most important of these arises directly from the definitions of the lower dimension and lower spectrum, which require a lower bound on the \textit{minimal} covering number. Such bounds are more challenging to obtain and can, in general, depend highly on the geometry of the underlying set. To overcome these issues, we start with an arbitrary value that satisfies the contrapositive assumption in the corresponding definitions in the source, and work towards such a value in the target. This facilitates distortion bounds by providing upper bounds on the corresponding minimal covering numbers in the target. However, it restricts the sequence of allowable scales from the start, which is not the case in \cite{OurQCspec}. The resulting lack of control over the gaps between consecutive levels is what leads to the different admissibility spectra exponents. 

A second (but less serious) complication comes from the lack of finite stability of the lower dimension and lower spectrum. The very first reduction in the proofs in \cite{OurQCspec} uses a Whitney decomposition and finite stability of $\dim_A$ and $\dim_{A,\reg}^\theta$ to focus attention on cubes that lie sufficiently far away from the boundary. This constraint allows for the application of the relevant reverse H\"older inequality. Since lower dimension and the lower spectrum are not finitely stable, this approach also must change. Instead, we consider appropriate thickenings of the source set and focus on cubes whose level is large enough to lie in such a thickening, far away from the boundary. 

Last but not least, in the case of the lower dimension and lower spectrum, we apply an improved reverse H\"older inequality a second time, relying on the self-improving Muckenhoupt property of the Jacobians of quasiconformal maps. This additional step results from the need for quantitative control of uniform constants between covering numbers and quotients of scales, which is not necessary in the case of the Assouad dimension and spectrum.

This paper is organized as follows. Section 2 contains the necessary background on dimension notions and quasiconformal mapping theory, as well as examples where the dimensions of different sets are evaluated. In Section 3 we prove Theorems \ref{thm:KI QC L dimension distortion} and \ref{thm:KI QC L spectrum distortion}, which imply Theorems \ref{thm: QC L dimension distortion} and \ref{thm: QC L spectrum distortion}, respectively. Section 4 contains the proofs of Theorems \ref{Thm: Hardy application R2} and \ref{th:ntalampekos-lower-bound}, along with other applications of Theorems \ref{thm: QC L dimension distortion} and \ref{thm: QC L spectrum distortion} on quasiconformal classification. In Section 5 we demonstrate how the same methods as in the proof of Theorems \ref{thm:KI QC L dimension distortion} yield corresponding estimates for the lower box dimension under Sobolev and quasiconformal maps.

\medskip

\paragraph{\bf Acknowledgements.} We thank Josh Kline and Ignacio Uriarte-Tuero for valuable conversations on the topics of this paper. We also thank Jonathan Fraser for helpful remarks on the general philosophy of dimension interpolation and specifically on the concepts of lower dimension and the lower spectrum. The first author is partially supported by an AMS-Simons Travel Grant. 

\medskip

\paragraph{\bf Use of Technology.} The authors used Wolfram Mathematica to expedite and verify calculations. Google search was also used, through which the `AI Overview' feature provided relevant information at times. Generative AI tools were not used at all in formulating or proving any of the results in this paper.

\section{Background}\label{sec:background}

For $x \in \R^d$ and $r>0$ we denote by $B(x,r)$ the (closed) Euclidean ball with center $x$ and radius $r$. Similarly, for $E\subset \R^d$ and $r>0$, we denote by $B(E,r)$ the closed $r$-neighborhood of $E$, i.e., the closure of the union $\bigcup_{x\in E}B(x,r)$. For a bounded set $E \subset \R^d$ and $r>0$ we denote by $N_r(E)$ the minimal number of sets of diameter at most $r$ needed to cover $E$.

For $x \in \R^d$ and $r>0$ we denote by $Q(x,r)$ the axes-parallel closed cube centered at $x$ with side length $2r$. In other words, if $x=(x_1,\ldots,x_d)$, then $Q(x,r) = [x_1-r,x_1+r] \times \cdots \times [x_d-r,x_d+r]$. Alternatively, $Q(x,r)$ is the closed metric ball in the $\ell^\infty$ norm on $\R^d$, with center $x$ and radius $r$.

For a set $E \subset \R^d$ and $x \in E$, $R>0$ we consider the axes-parallel cube $Q:=Q(x,R) \subset \R^d$. We subdivide $Q$ into $2^d$ essentially disjoint subcubes, each with side length equal to half of the side length of $Q$, and then we subdivide each of those cubes in the same fashion, and so on. Let $\cW(Q)$ denote the collection of all such cubes obtained at any level of the construction, and let $\cW_m(Q)$ denote the collection of all cubes obtained after $m$ steps. The collection $\cW(Q)$ is known as the {\it dyadic decomposition} of $Q$, while $\cW_m(Q)$ consists of the {\it level $m$ dyadic subcubes} in $Q$. 

We will make use of dyadic versions of the covering numbers $N_r(E)$ defined above. For a set $E \subset \R^d$ and for $m \in \N$ denote by $N(B(x,R) \cap E,m)$ the number of level $m$ dyadic subcubes in $Q$ needed to cover $B(x,R) \cap E$. We record the following comparability statements relating the usual covering numbers and their dyadic counterparts; these will be used in the proof of Proposition \ref{prop:lower-spectrum-facts}. For $x \in E$, $R>0$, and $m \in \N$, we have
	\begin{equation}\label{eq: cov numbers ineq}
		N_{2^{-m+1}R \sqrt{d}}(B(x,R)\cap E) \leq N(B(x,R)\cap E,m) \leq 3^d N_{2^{-m+1}R}(B(x,R)\cap E).
	\end{equation}
    To see why \eqref{eq: cov numbers ineq} holds true, note first that since cubes of side length $2^{-m+1}R$ have diameter at most $2^{-m+1}R \sqrt{d}$, the left hand inequality is trivial. For the right hand inequality it suffices to note that every set of diameter at most $2^{-m+1}R$ cannot possibly intersect more than $3^d$ axes-parallel cubes of side length $2^{-m+1}R$.

Finally, we record the classical {\it Morrey--Sobolev inequality} valid for continuous functions $f:\Omega \to \R^d$ defined on a domain $\Omega \subset \R^d$ and lying in the Sobolev class $W^{1,p}_\loc(\Omega:\R^d)$ for some $p>d$: there exists a constant $C_{MS} = C_{MS}(d,p) > 0$ so that for any closed cube $Q(x,R)$ contained in $\Omega$ and all $y,z \in Q(x,R)$ the estimate
\begin{equation}\label{eq:morreysobolev}
|f(y)-f(z)| \le C_{MS} |y-z|^{1-d/p} \left( \int_{Q(x,R)} |Df|^p \right)^{1/p}
\end{equation}
holds true. Here, for a mapping $f:\Omega \to \R^d$ lying in the local Sobolev space $W^{1,1}_\loc(\Omega:\R^d)$, we denote by $Df$ the (a.e.\ defined) differential matrix and by $|\bA|=\max\{|\bA(\bv)|:|\bv|=1\}$ the operator norm of a matrix $\bA$.

\subsection{Quasiconformal mappings.}\label{subsec:qc-mappings}

A homeomorphism $f:\Omega \to \Omega'$ between domains in $\R^d$, $d \ge 2$, is {\it $K-$quasiconformal}, $K\geq1$, if $f$ lies in the local Sobolev space $W^{1,d}_\loc(\Omega:\R^d)$ and satisfies the distortion inequality
\begin{equation}\label{eq:qc-defn}
|Df|^d \le K \det Df
\end{equation}
a.e.\ in $\Omega$. The smallest value $K \ge 1$ for which \eqref{eq:qc-defn} holds a.e.\ in $\Omega$ is known as the {\it outer dilatation} of $f$ and is denoted $K_O(f)$. Alternatively, for an invertible matrix $\bA$ set $\ell(\bA):=\min\{|\bA(\bv)|:|\bv|=1\}$ and replace \eqref{eq:qc-defn} with the inequality
\begin{equation}\label{eq:qc-defn-2}
\det Df \le K \ell(Df)^d \,.
\end{equation}
The smallest choice of $K$ for which \eqref{eq:qc-defn-2} holds a.e.\ is the {\it inner dilatation} of $f$ and is denoted $K_I(f)$. We set $K_f := \max\{K_O(f),K_I(f)\}$. The outer and inner dilatations are related by the mutual inequalities $K_O(f) \le K_I(f)^{d-1}$ and $K_I(f) \le K_O(f)^{d-1}$; thus \begin{equation}\label{eq:dimension-2-Ks}
K_O(f)=K_I(f):=K_f
\end{equation}
when $d=2$.

Henceforth we denote by $\Barint_E w$ the average of a locally integrable function $w:\R^d \to \R$ over a set $E$ of positive measure, namely, $\Barint_E w = |E|^{-1} \int_E w$. We next recall the reverse H\"older inequality for the differential matrix of a quasiconformal mapping, first established by Gehring in \cite{Gehring73}. Let $f:\Omega \to \Omega'$ be a $K$-quasiconformal mapping between domains in Euclidean space $\R^d$, $d \ge 2$. According to Gehring's well-known higher integrability theorem \cite[Theorem 1]{Gehring73}, $f$ lies in a local Sobolev space $W^{1,p}_\loc(\Omega:\R^d)$ for some $p>d$. In particular, the differential of $f$ satisfies a reverse H\"older inequality of the form 
\begin{equation}\label{eq:pnK0}
\left( \Barint_Q |Df|^p \right)^{1/p} \le C_{\RH} \left( \Barint_{Q} |Df|^d \right)^{1/d} 
 \end{equation}
for all cubes $Q \subset \Omega$ with $\diam(Q) < \dist(Q,\partial\Omega)$ and $\diam(f(Q)) < \dist(f(Q),\partial\Omega')$. Here $C_{\RH}>0$ denotes a constant which depends only on the dimension $d$ and the dilatation $K$ of $f$.  

For $d \ge 2$ and $K \ge 1$, denote by
$$
p_O^{\RH}(d,K)
$$
the supremum of those values $p>d$ so that \eqref{eq:pnK0} holds true for every quasiconformal map $f$ between domains in $\R^d$ with $K_O(f) \le K$. The superscript in the notation for this higher integrability exponent indicates the role of the reverse H\"older inequality in the definition. We also denote by
$$
p_O^{\Sob}(d,K)
$$
the supremum of those values $p>d$ so that every quasiconformal map $f$ between domains $\Omega$ and $\Omega'$ in $\R^d$ with $K_O(f) \le K$ necessarily lies in $W^{1,p}_\loc(\Omega:\R^d)$. The value $p_O^{\Sob}(d,K)$ is what is traditionally referred to as the higher integrability exponent for quasiconformal mappings in the literature, where it is generally denoted $p(d,K)$. Note that
\begin{equation}\label{eq:pRHSobnK}
p_O^{\RH}(d,K) \le p_O^{\Sob}(d,K)
\end{equation}
and moreover,
\begin{equation}\label{eq:pRHSobnK2}
p_O^{\Sob}(d,K) \le \frac{dK}{K-1}
\end{equation}
as witnessed by the radial stretch map $f(x) = |x|^{1/K-1}x$. In contrast with the situation for Hausdorff or upper box dimension, distortion estimates for the Assouad or lower dimension are phrased in terms of the reverse H\"older inequality higher integrability exponent $p_O^{\RH}(d,K)$. However, it is widely expected that equality holds in both \eqref{eq:pRHSobnK} and \eqref{eq:pRHSobnK2} for all $d \ge 3$ and $K \ge 1$.

In a similar fashion, let $p_I^{\RH}(d,K)$, respectively $p_I^{\Sob}(d,K)$, denote the supremum of values so that \eqref{eq:pnK0} holds true, respectively $f$ lies in $W^{1,p}_\loc(\Omega:\R^d)$, for every quasiconformal map $f$ between domains $\Omega$ and $\Omega'$ in $\R^d$ with $K_I(f) \le K$. When $d=2$ the identities in \eqref{eq:dimension-2-Ks} imply that the higher integrability exponents $p^{\RH}(2,K)$ and $p^{\Sob}(2,K)$ are defined without reference to either outer or inner dilatation. Astala \cite{Astala} showed that
\begin{equation}\label{eq:p2K}
p^{\Sob}(2,K) = \frac{2K}{K-1}
\end{equation}
and the stronger statement
\begin{equation}\label{eq:p2K2}
p^{\RH}(2,K) = \frac{2K}{K-1}
\end{equation}
can be found in \cite[Theorem 13.2.3 and Corollary 13.2.4]{aim:2qc}.

The following property of quasiconformal mappings is popularly known as the `egg yolk principle'. See, for instance, \cite[Theorem 11.14]{hei:lectures}.
\begin{proposition}\label{prop:egg-yolk}
Fix $d \ge 2$ and $K \ge 1$. Then there exists an increasing homeomorphism $\eta = \eta_{K,d}$ of $[0,\infty)$ so that for any $K$-quasiconformal homeomorphism $f:\Omega \to \Omega'$ between domains in $\R^d$ and any cube $Q \subset \Omega$ with $\diam(Q) \le \dist(Q,\partial\Omega)$, the restriction $f|_Q$ is $\eta$-quasisymmetric. More precisely, if $x,y,z \in Q$ and $x \ne z$, then
$$
\frac{|f(x)-f(y)|}{|f(x)-f(z)|} \le \eta \left( \frac{|x-y|}{|x-z|} \right) \, .
$$
\end{proposition} In later arguments, it is the above property that allows us to relate radii of balls in the source with diameters of images of said balls under quasiconformal maps.

\subsection{Assouad and lower dimensions}

We first recall the definitions of Assouad and lower dimension for sets in Euclidean space. A good reference for this material is \cite[Chapters 2 and 3]{FraserBook}.

Let $E$ be a non-empty subset of $\R^d$. The {\it Assouad dimension} of $E$, denoted $\dim_A E$, is the infimum of all values $s>0$ for which there exists a constant $C>0$ so that
$$
N_r(B(x,R)\cap E) \le C \left( \frac{R}{r} \right)^s
$$
for all $x \in E$ and all $0 < r \le R$. Analogously, the {\it lower dimension} of $E$, denoted $\dim_L E$, is the supremum of all values $s \ge 0$ for which there exists a constant $C>0$ so that 
$$
N_r(B(x,R)\cap E) \ge (1/C) \left( \frac{R}{r} \right)^s
$$
for all $x \in E$ and all $0 < r \le R<\diam E$. We point out that some references call these two concepts the {\it upper Assouad dimension} and {\it lower Assouad dimension} respectively.

The Assouad dimension of a set $E$ captures the location $x$ and scale $R$ in $E$ where the set is `thickest' in the sense that the relative covering number $N_r(B(x,R)\cap E)$ exhibits the largest possible growth as a function of $\tfrac{R}r$. On the other hand, the lower dimension of $E$ captures the location and scale where the set is `thinnest' in the sense that such covering number exhibits the smallest possible growth in terms of $\tfrac{R}r$. 

The following proposition records some basic observations about these two notions of dimension. For more details, see \cite[Chapters 2 and 3]{FraserBook}.

\begin{proposition}\label{prop:Assouad_and_lower_properties}
\begin{itemize}
    \item[(a)] The Assouad and lower dimensions are both stable under closure: if $E \subset \R^d$ then $\dim_A E = \dim_A \overline{E}$ and $\dim_L E = \dim_L \overline{E}$.
    \item[(b)] Both notions of dimension are invariant under bi-Lipschitz mappings. However, each of these notions can in general increase under Lipschitz mappings.
    \item[(c)] The Assouad dimension fails to be countably stable, but it is finitely stable:
    $$
    \dim_A(E \cup F) = \max\{\dim_A E,\dim_A F\}
    $$
    for all $E,F \subset \R^d$. The lower dimension satisfies the finite stability condition
    $$
    \dim_L(E \cup F) = \min\{\dim_L E,\dim_L F\}
    $$
    for all $E,F \subset \R^d$ such that $\dist(E,F)>0$. However, for arbitrary sets $E$ and $F$ there is in general no relationship between $\dim_L(E\cup F)$ and the lower dimensions of $E$ and $F$.
    \item[(d)] The Assouad dimension is monotonic: if $E \subset F$ then $\dim_A E \le \dim_A F$. The corresponding statement for the lower dimension need not hold true in general.
    \item[(e)] A set $E\subset\R^d$ is porous if and only if $\dim_A E < d$. Moreover, $E$ is uniformly perfect if and only if $\dim_L E > 0$.
    \item[(f)] The Assouad dimension satisfies the {\it open set property}: if $U$ is a nonempty open subset of $\R^d$ then $\dim_A U = d$. There exist open sets $U \subset \R^d$ for which $\dim_L U < d$.
\end{itemize} 
\end{proposition}

The failure of monotonicity for lower dimension in part (d) may seem strange to the uninitiated reader. As an example, we invite the reader to compute the lower dimension of the subsets $E = [0,1]$ and $F = E \cup \{2\}$ of $\R$.

Recall that a set $E \subset \R^d$ is {\it porous} if there exists a constant $c>0$ so that for every $x \in E$ and $0<R<\diam E$ there exists $z \in B(x,R)$ so that $B(z,cR) \cap E = \emptyset$. On the other hand, $E$ is {\it uniformly perfect} if there exists $c>0$ so that for every $x \in E$ and $0<R<\diam E$, the set $B(x,R) \cap E \setminus B(x,cR)$ is nonempty. The conclusions in part (e) were proved by Luukkainen \cite{Luukkainen} (for the Assouad dimension) and by K\"aenm\"aki, Lehrb\"ack and Vuorinen \cite{klv:whitney} (for the lower dimension). See also \cite[Chapter 13]{FraserBook} for more information.

The claim about lower dimensions in part (f) does not appear explicitly in \cite{FraserBook}. We provide examples below in section \ref{subsec:examples}.

Parts (c), (d) and (f) of the above proposition indicate that the lower dimension fails to enjoy many natural properties typically associated with a concept of dimension. Nevertheless, the lower dimension plays an important role in dimension theory; see, for instance, \cite[Theorem 1.3]{ge:trees} for an application to embedding theory or \cite[Proposition 3.4]{bhm:quasivisual} for an application to uniformization and dynamics. In this paper we show that it also transforms naturally under quasiconformal mappings.

Porosity \cite[Theorem 3.1]{vai:porosity} and uniform perfectness \cite[Corollary 4.6]{UnifPerfectInv} are invariant under global quasiconformal mappings. In view of Proposition \ref{prop:Assouad_and_lower_properties} (e), the strict inequalities $\dim_A E < d$ and $\dim_L E > 0$ are preserved by global quasiconformal mappings. The results in \cite{OurQCspec} together with Theorem \ref{thm: QC L dimension distortion} provide quantitative strengthenings of this conclusion for the two notions of dimension.

\subsection{The lower spectrum}

The Assouad spectrum and lower spectrum, introduced in \cite{fy:assouad-spectrum}, are two examples of the general framework of dimension interpolation pioneered by Fraser, Falconer, and their collaborators. In this paper we are particularly interested in the lower spectrum, whose definition we now recall. Let $E$ be a non-empty subset of $\R^d$. For $0<\theta<1$, define
$$
\dim_L^\theta E
=
\sup\left\{
\alpha\geq 0:
\begin{array}{l}
\exists \, C>0 \text{ such that for all } x\in E
\text{ and all } 0<r<1,\\[2mm]
N_r\bigl(B(x,r^\theta)\cap E\bigr)
\geq
(1/C) \left(\dfrac{r^\theta}{r}\right)^\alpha
\end{array}
\right\}.
$$
The function  $\dim_L^\theta E:(0,1)\to [0,d]$ comprises the lower spectrum of $E$.

In studying dimension distortion properties under Sobolev and quasiconformal mappings, certain modified versions of these dimension notions are more convenient. This principle was already articulated in \cite{OurQCspec} in the context of quasiconformal distortion of Assouad-type dimensions.

For $0<\theta<1$, define
$$
\dim_{L,\text{reg}}^\theta E
=
\sup\left\{
\alpha\geq 0:
\begin{array}{l}
\exists\, C>0 \text{ such that for all } x\in E
\text{ and all } 0<r\leq R^{1/\theta}<R<1,\\
N_r\bigl(B(x,R)\cap E\bigr)
\geq
C\left(\dfrac{R}{r}\right)^\alpha
\end{array}
\right\}.
$$
We call $\dim_{L,\text{reg}}^\theta E$ the {\it regularized lower spectrum of $E$ with parameter $\theta$} to match the terminology from \cite{OurQCspec}. However, the term \textit{quasi-lower} spectrum has also been used in the literature \cite{ChenLowerSpec}.

The regularized spectra can be realized as monotonic envelopes of their unregularized counterparts. For the Assouad dimension case see \cite[Theorem 3.3.6]{FraserBook}. The following proposition has been proved by Hare and Troscheit \cite[Theorem A.2]{KatTros_lower_dim_spec} in the general context of doubling metric spaces. For the convenience of the reader, we include a simplified version of the proof in the Euclidean setting.

\begin{proposition}
Let $E\subset \R^d$. Then, for all $\theta\in (0,1)$ we have
    $$
\dim_{L,\text{reg}}^\theta E
=
\inf_{0<\phi\leq\theta} \dim_L^\phi E.
$$
\end{proposition}

In particular, the function
$$
(0,1) \ni \theta \mapsto \dim_{L,\text{reg}}^\theta E
$$
is monotonically {\bf non-increasing}. (Recall that the regularized Assouad spectrum $(0,1) \ni \theta \mapsto \dim_{A,\text{reg}}^\theta E$ is monotonically non-decreasing.)

\begin{proof}
    Fix $\theta \in (0,1)$ and let $D$ denote the infimum on the right hand side of the equality.
    Set
    $$
    A_\theta=\left\{
    \alpha\geq 0:
    \begin{array}{l}
    \exists\, C>0 \text{ such that for all } x\in E
    \text{ and all } 0<r\leq R^{1/\theta}<R<1,\\
    N_r\bigl(B(x,R)\cap E\bigr)
    \geq
    (1/C)\left(\dfrac{R}{r}\right)^\alpha
    \end{array}
    \right\},
    $$
    and, for every $\phi\in (0,\theta)$ set
    $$
    B_\phi=\left\{
    \alpha\geq 0:
    \begin{array}{l}
    \exists \, C>0 \text{ such that for all } x\in E
    \text{ and all } 0<r<1,\\[2mm]
    N_r\bigl(B(x,r^\phi)\cap E\bigr)
    \geq
    (1/C) \left(\dfrac{r^\phi}{r}\right)^\alpha
    \end{array}
    \right\}.
    $$
    Let $s\in A_\theta$, $x\in E$, $r\in (0,1)$, and $\phi\in (0,\theta)$. Then for $R=r^\phi<1$ we have
    $0<r\leq r^{\phi/\theta}=R^{1/\theta}<R<1$, and so there is $C>0$ such that
    $$
    N_r(B(x,R)\cap E)\geq (1/C) \left( \frac{R}{r}\right)^s.
    $$
    Since $R=r^\phi$, this implies that $s\in B_\phi$. Hence $A_\theta\subset B_\phi$ for any $\phi\in (0,\theta)$, which yields
    \begin{equation}\label{eq:an-upper-bound}
    \dim_{L,\text{reg}}^\theta E
    \leq D.
    \end{equation}
    Observe that equality in \eqref{eq:an-upper-bound} holds trivially if $D=0$.  
    For the remainder of the proof we assume that $D>0$.
    Let $s\in (0, D)$
    and fix $u\in (s, D)$.
    Set
    $$
    \phi_j:= 1-\left(\frac{u}{s}\right)^j (1-\theta),
    $$ for all $j\in \{0, 1, \dots, N\}$, where $N$ is the largest integer for which $\phi_j>0$. Note that for all $j$ we have
    \begin{equation}\label{eq: phi_j+1}
        \phi_{j+1}=1-\frac{u}{s} (1-\phi_j),
    \end{equation} and also $\phi_j\in (0,\theta)$.

    Since $u<\dim_L^\phi E$ for all $\phi\in(0,\theta)$, we have that there exist $C_0, C_1, \dots, C_N>0$ such that
    \begin{equation}\label{eq: Nr for phi_j}
        N_r(B(x,r^{\phi_j})\cap E)\geq (1/C_j) \left(\frac{r^{\phi_j}}{r}\right)^u,
    \end{equation}
    for all $j\in \{0,\dots, N\}$ and all $x\in E$, $r\in (0,1)$. Set $C=\max\{C_0, \dots, C_N\}$.

    Let $x\in E$ and $0<r\leq R^{1/\theta}<R<1$. Set
    $$
    \phi=\frac{\log R}{\log r}\in (0,1).
    $$ 
    Then $R=r^\phi$ and $\phi\in (0,\theta]$. Fix $j\in \{0,\dots, N\}$ such that
    \begin{equation}\label{eq:phi-j-phi}
    \phi_{j+1}\leq \phi<\phi_j.
    \end{equation} 
    Then \eqref{eq:phi-j-phi} and \eqref{eq: phi_j+1} imply
    \begin{equation}\label{eq: s and u exponents}
        s(1-\phi)\leq s(1-\phi_{j+1})= u(1-\phi_j)
    \end{equation}
    Since $r^{\phi_j}\leq r^\phi=R$, we have that
    $$
    N_r(B(x,R)\cap E)\geq N_r(B(x,r^{\phi_j})\cap E).
    $$
    By \eqref{eq: Nr for phi_j} and the choice of $C$ implies
    $$
    N_r(B(x,R)\cap E)\geq (1/C) \left(\frac{r^{\phi_j}}{r}\right)^u.
    $$ 
    But the right-hand side is equal to $(1/C)(1/r)^{u(1-\phi_j)}$, and so the above and \eqref{eq: s and u exponents} yield
    $$
    N_r(B(x,R)\cap E)\geq (1/C) \left(\frac{1}{r}\right)^{s(1-\phi)}=(1/C) \left(\frac{R}{r}\right)^{s}.
    $$ 
    Since $x$, $r$, and $R$ are arbitrary and $C$ is independent of the choice of the point and scales, the above implies that $s<\dim_{L,\text{reg}}^\theta E$. Letting $s \nearrow D$
    completes the proof.
\end{proof}

The lower dimension and lower box-counting dimensions bound the (regularized) lower spectrum as follows:
$$
\dim_L E \le \dim_{L,reg}^\theta E \le \underline{\dim}_B E, \qquad \forall \, 0 < \theta < 1.
$$
In particular, $\lim_{\theta \to 0} \dim_{L,reg}^\theta E \le \underline{\dim}_B E$ and $\lim_{\theta \to 1} \dim_{L,reg}^\theta E \ge \dim_L E$. However, both of the latter two inequalities can be strict. For examples, we refer to \cite{fy:assouad-spectrum} and \cite{FraserBook}.

In mapping-theoretic applications it is important to allow some flexibility in the radius bound $r \le R^{1/\theta}$ due for instance to the presence of multiplicative constants in the H\"older continuity estimates for quasiconformal mappings. Part (i) of the following proposition guarantees that such modifications do not alter the lower spectrum. In part (ii) of the proposition, we characterize the lower spectrum in terms of the growth of the dyadic covering numbers $N(B(x,R)\cap E,m)$.

\begin{proposition}\label{prop:lower-spectrum-facts}
Let $E \subset \R^d$ and $0<\theta<1$.
\begin{itemize}
\item[(i)] Fix a value $C_1 \ge 1$. Then the regularized lower spectrum value $\dim_{L,\text{reg}}^\theta E$ is equal to the supremum of all $\alpha\geq 0$ for which there exists $C>0$ so that
    $$
    N_r(B(x,R) \cap E) \geq (1/C)\left(\frac{R}{r}\right)^\alpha,
    $$
    for all $x \in E$ and all $0<r/C_1\le R^{1/\theta} < R < 1$.

\item[(ii)] The regularized lower spectrum value $\dim_{L,\text{reg}}^\theta E$ is equal to the supremum of all $\alpha\geq 0$ for which there exists $C>0$ so that
    $$
    N(B(x,R) \cap E,m) \ge (1/C) 2^{m\alpha},
    $$
    for all $x \in E$ and all $m\in \N, R>0$ with $0<2^{-m}R\le R^{1/\theta} < R < 1$.
\end{itemize}
\end{proposition}

\begin{proof}
	(i) Denote by 
    $$
    A_\theta=\left\{
    \alpha\geq 0:
    \begin{array}{l}
    \exists\, C>0 \text{ such that for all } x\in E
    \text{ and all } 0<r\leq R^{1/\theta}<R<1,\\
    N_r\bigl(B(x,R)\cap E\bigr)
    \geq
    (1/C)\left(\dfrac{R}{r}\right)^\alpha
    \end{array}
    \right\},
    $$
	and by
	$$
    B_\theta=\left\{
    \alpha\geq 0:
    \begin{array}{l}
    \exists\, C>0 \text{ such that for all } x\in E
    \text{ and all } 0<r/C_1\leq R^{1/\theta}<R<1,\\
    N_r\bigl(B(x,R)\cap E\bigr)
    \geq
    (1/C) \left(\dfrac{R}{r}\right)^\alpha
    \end{array}
    \right\},
    $$
	By the definition of the regularized lower spectrum $\dim_{L,\text{reg}}^\theta E=\sup A_\theta$, so it is enough to show that $A_\theta=B_\theta$.
	
	Let $\alpha\in A_\theta$, $x \in E$ and $r, R>0$ with $0<r/C_1\le R^{1/\theta}<R<1$. It is a standard
    property of $\R^d$ that there exists a constant $M=M(d,C_1)\geq 1$ such that every subset of $\R^d$ of diameter at most $r$ can be covered by at most $M$ sets of diameter at most $r/C_1$. Therefore, 
    $$
	N_{r/C_1}(B(x,R)\cap E) \leq M N_r(B(x,R)\cap E),
	$$ which by definition of $A_\theta$ implies
	$$
    N_r(B(x,R)\cap E) \geq  (CM)^{-1} \left(\frac{R}{r/C_1}\right)^{\alpha} = (CM)^{-1} C_1^\alpha (R/r)^{\alpha}.
    $$
	Hence $\alpha\in B_\theta$ and so $A_\theta \subset B_\theta$.
	
	Since $C_1\geq 1$, for any $0<r\leq R^{1/\theta} < R<1$ we have that $r/C_1 \le r\leq R^{1/\theta}$, which makes the inclusion $B_\theta \subset A_\theta$ trivial. Therefore, $A_\theta = B_\theta$ and the proof of (i) is complete.
	
	(ii) Denote by 
    $$
    D_\theta := \left\{
    \alpha \geq 0 \,:\,
    \begin{array}{l}
    \exists\, C>0 \text{ such that }
    N(B(x,R) \cap E,m) \geq (1/C) 2^{m\alpha}\\
    \text{for all } m\in\N,\ R>0 \text{ with }
    0<2^{-m}R\leq R^{1/\theta}<R<1\,
    \text{and all } x\in E
    \end{array}
    \right\}.
    $$
	
	In a similar fashion as in previous arguments, it suffices to show that $A_\theta = D_\theta$. 
	
	Now let $\alpha\in A_\theta$.  Let $x \in E$ and $m \in \N$,  $R>0$ with $0<2^{-m}R\le R^{1/\theta} < R < 1$. Then
	$$
	\frac{2^{-m+1}R\sqrt{d}}{2\sqrt{d}}=2^{-m}R\leq R^{1/\theta},
	$$
    which by part (i) for  $C_1=2\sqrt{d}$ implies
	$$
	N_{2^{-m+1}R\sqrt{d}}(B(x,R)\cap E) \geq (1/C) \left( \frac{R}{2^{-m+1}R\sqrt{d}}\right)^\alpha
	= (1/C)(2\sqrt{d})^{-\alpha}2^{m\alpha}.
	$$
	By \eqref{eq: cov numbers ineq}, this implies that
	$$
	N(B(x,R)\cap E,m) \geq (1/C) (2\sqrt{d})^{-\alpha}2^{m\alpha}.
	$$
	Hence $\alpha \in D_\theta$ and so $A_\theta \subset D_\theta$.
	
	Finally, let $\alpha \in D_\theta$, $x \in E$ and $r, R>0$ with $0<r\leq R^{1/\theta}<R<1$. Let $m \in \N$ be such that
	$$
	2^{-m}R\leq r<2^{-m+1}R.
	$$
	Then $0<2^{-m}R\leq R^{1/\theta}<R<1$, and due to $\alpha\in D_\theta$ we have
	$$
	N(B(x,R)\cap E,m) \geq (1/C) 2^{m\alpha}.
	$$
	Moreover, by \eqref{eq: cov numbers ineq} we have
	$$
	N_{2^{-m+1}R}(B(x,R)\cap E) \geq 3^{-d}N(B(x,R)\cap E,m),
	$$
	and since $r<2^{-m+1}R$, it follows that
	$$
	N_r(B(x,R) \cap E) \geq N(B(x,R) \cap E,2^{-m+1}R) \geq 3^{-d} (1/C) 2^{m\alpha}.
	$$
	Since $2^{-m}R\leq r$ the above implies
	$$
	N_r(B(x,R) \cap E) \geq 3^{-d}(1/C)(R/r)^\alpha.
	$$
	Hence $\alpha \in A_\theta$ and so $D_\theta \subset A_\theta$. This completes the proof.	
\end{proof}

\subsection{Examples}\label{subsec:examples}

The fact that the lower dimension (and lower spectrum) of an open subset of $\R^d$ can be strictly less than $d$ is disconcerting, when viewed from the perspective of more classical dimensional notions such as Hausdorff or Assouad dimension. However, we obtain interesting classification results by coupling this fact with our quasiconformal distortion theorems. In this section, we compute the lower dimension and lower spectrum of a variety of open sets.

\begin{example}\label{ex:fraser-yu}
Our first example appears as Example 2.5 in \cite{Fraser:lower-dimension}. The open set $V = \cup_{n=1}^\infty (\tfrac1n - 2^{-n},\tfrac1n + 2^{-n}) \subset \R$ has lower dimension equal to zero. In fact, the proof of this fact in \cite[Example 2.5]{Fraser:lower-dimension} shows that $\dim_{L,reg}^\theta V = 0$ for every $\theta\in(0,1)$. Note, however, that the lower box dimension of $V$ is equal to $1$.
\end{example}

By \cite[Theorem 2.1]{Fraser:lower-dimension} and \cite[Proposition 4.5]{fy:assouad-spectrum}, it follows that the $d$-fold product $V^d$ is a nonempty open subset of $\R^d$ with lower dimension and lower spectrum identically equal to zero. By \cite[Theorem 2.1]{Fraser:lower-dimension} and \cite[Proposition 4.4]{fy:assouad-spectrum}, we also have that $V^{d-k} \times \R^k$ is a nonempty open subset of $\R^d$ with lower dimension and lower spectrum identically equal to $k$. Here $k$ is any integer in $\{1,\ldots,d-1\}$.

We next give examples of planar open sets with nonconstant lower spectrum.

\begin{example}\label{ex:cusp-domains}
Fix $\alpha>1$ and set $\Omega_\alpha := \{z = x + \bi y : 0 <x<1, |y| < x^\alpha\}$. We claim that $\dim_L(\Omega_\alpha) = 1$ and
$$
\dim_{L,reg}^\theta(\Omega_\alpha) = d_{\alpha,\theta}
$$
where
$$
d_{\alpha,\theta} := \begin{cases} 
\frac{2-(1+\alpha)\theta}{1-\theta}, & 0<\theta<\tfrac1\alpha, \\ 1, & \tfrac1\alpha \le \theta < 1. \end{cases}
$$
First, we show that 
\begin{equation}\label{eq:dimLregOmegaalphaupperbound}
\dim_{L,reg}^\theta(\Omega_\alpha) \le  d_{\alpha,\theta}.
\end{equation}
Let $t \in (0, \dim_{L,reg}^\theta(\Omega_\alpha))$. Then there exists $C>0$ so that for every $z \in \Omega_\alpha$ and every $0<r\le R^{1/\theta} < R \le 1$ it holds true that
\begin{equation}\label{eq:here-is-a-lower-bound}
N_r(B(z,R)\cap \Omega_\alpha) \ge \frac1C (\frac{R}{r})^t.
\end{equation}
We choose $z=0$ and consider the sequences $R_n = \tfrac1n$ and $r_n = \tfrac1{n^{1/\theta}}$. Observe that $N_{r_n}(B(z,R_n)\cap \Omega_\alpha) \le N_{r_n}(W)$ where $W = \{z = x + \bi y \in \Omega:0<x<R_n\}$. We estimate $N_{r_n}(W)$ from above by considering an essentially disjoint collection of closed squares of side length $r_n/\sqrt2$ covering $W$. Note that each such square is contained within a ball of radius $r_n$, so estimating the cardinality of such a collection of squares from above gives an upper bound for $N_{r_n}(W)$. The minimal collection of such squares drawn from the standard rectilinear tiling of the plane by squares of side length $r_n/\sqrt2$ has cardinality
$$
\lesssim \sum_{j=1}^M \left( \frac{(jr_n)^\alpha}{r_n} + 1 \right),
$$
where $M$ is an integer satisfying $M \simeq R_n/r_n$. In fact, it suffices to consider an upper Riemann sum for the integral $\int_0^{R_n} 2x^\alpha \, dx$ built from squares in the previously mentioned rectilinear tiling. Using the values of $R_n$ and $r_n$ we conclude that
\begin{equation}\label{eq:here-is-an-upper-bound}
N_{r_n}(B(z,R_n)\cap \Omega_\alpha) \le N_{r_n}(W) \lesssim \frac{R_n^{1+\alpha}}{r_n^2} + \frac{R_n}{r_n} = n^{-1-\alpha+\tfrac2\theta} + n^{-1+\tfrac1\theta}.
\end{equation}
Combining \eqref{eq:here-is-a-lower-bound} and \eqref{eq:here-is-an-upper-bound} yields
$$
n^{(\tfrac1\theta-1)t} \lesssim n^{-1-\alpha+\tfrac2\theta} + n^{-1+\tfrac1\theta}.
$$
If $\theta>\tfrac1\alpha$ then the dominant term on the right hand side is $n^{-1+\tfrac1\theta}$ and the inequality can hold for all $n$ only if $t \le 1$. On the other hand, if $\theta<\tfrac1\alpha$ then the dominant term is $n^{-1-\alpha+\tfrac2\theta}$ and the inequality can hold for all $n$ only if $(\tfrac1\theta-1)t \le -1-\alpha+\tfrac2\theta$. These conclusions imply that we must have $t \le d_{\alpha,\theta}$. Letting $t \nearrow \dim_{L,reg}^\theta(\Omega_\alpha)$ we obtain \eqref{eq:dimLregOmegaalphaupperbound}.

Next, we show that 
\begin{equation}\label{eq:dimLregOmegaalphalowerbound}
\dim_{L,reg}^\theta(\Omega_\alpha) \ge  d_{\alpha,\theta}.
\end{equation}
Towards a contradiction, assume that $d_{\alpha,\theta}>\dim_{L,reg}^\theta(\Omega_\alpha)$. Then there exist points $z_n \in \Omega_\alpha$ and scales $0<r_n \le R_n^{1/\theta} < R_n \le 1$ so that
$$
N_{r_n}(B(z_n,R_n)\cap \Omega_\alpha) \le \frac1n \, \left( \frac{R_n}{r_n} \right)^{d_{\alpha,\theta}}.
$$
By studying the area of $B(z,R_n)\cap\Omega_\alpha$ as a function of $z$ we note that it suffices to assume that $z=0$. We now estimate the quantity $N_{r_n}(B(0,R_n) \cap \Omega_\alpha)$ from below. Note that the cardinality of any covering of $B(0,R_n) \cap \Omega_\alpha$ by balls of radius $r_n$ is bounded below by the cardinality of any $3r_n$-net of points in $B(0,R_n) \cap \Omega_\alpha$. Thus
$$
N_{r_n}(B(0,R_n) \cap \Omega_\alpha) \ge M_n,
$$
where $M_n$ is the maximal cardinality of such a $3r_n$-net. Consider the points $q_{j,\ell} = 3jr_n + 3 \ell r_n \bi$, where  $j=0,\ldots,M$ and $\ell = -L_j,\ldots,L_j$. Here $M$ is an integer such that $M \simeq R_n/r_n$, while $L_j$ is an integer such that $L_j \simeq j^\alpha \, r_n^{\alpha-1}$. The total number of points is
$$
\sum_{j=0}^M (2L_j+1) = M + 2 \sum_{j=1}^M L_j \simeq \frac{R_n}{r_n} + \frac{R_n^{1+\alpha}}{r_n^2}  = n^{-1+\tfrac1\theta} + n^{-1-\alpha+\tfrac2\alpha}.
$$
We conclude that
$$
n^{-1+\tfrac1\theta} + n^{-1-\alpha+\tfrac2\alpha} \lesssim \frac1n \, n^{(\tfrac1\theta-1)d_{\alpha,\theta}}
$$
and we again reach a contradiction by considering the expression for $d_{\alpha,\theta}$.
\end{example}

\begin{example}\label{ex:another-example}
Fix $0<b<a<\tfrac12$ and let $(\beta_j)$ and $(\varphi_j)$ be sequences satisfying 
\begin{equation}\label{eq:beta-j-bound}
0<\beta_j \le C_0 b^j
\end{equation}
and 
\begin{equation}\label{eq:varphi-j-bound}
0<\varphi_j \le C_0 a^{-j} \beta_j
\end{equation}
for some $C_0$. 

For $j \in \N$ denote by $S_j$ the sector consisting of points $z \in \C$ such that $a^j < |z| < a^j + \beta_j$ and $|\arg(z)|< \varphi_j$. We consider only those values $j \ge N$ where $N \in \N$ satisfies
\begin{equation}\label{eq:N}
(\frac{b}{a})^N < \frac{\tfrac12-a}{C_0b}.
\end{equation}

\begin{lemma}
For $j \ge N$ with $N$ chosen as in \eqref{eq:N} it holds true that
\begin{equation}\label{eq:distance-lower-bound}
\dist(S_j,S_{j+1}) \ge \frac12 a^j
\end{equation}
and
\begin{equation}\label{eq:diameter-upper-bound}
\diam(S_j) \le C_* \beta_j.
\end{equation}
Here $C_*$ denotes a constant which depends only on $C_0$.
\end{lemma}

\begin{proof}
Let $j \ge N$. If $z \in S_j$ and $z' \in S_j$, then
$$
|z-z'| \ge \biggl| |z| - |z'| \biggr| \ge a^j - (a^{j+1} - \beta_{j+1}) \ge (1-a) a^j - C_0 b^{j+1} \ge \frac12 a^j
$$
by the choice of $N$. If $z,\zeta \in S_j$ then
\begin{equation*}\begin{split}
|z-\zeta|^2 
&= (|z|-|\zeta|)^2 + 2 |z|\,|\zeta|(1-\cos(\arg(\overline{z}\zeta))) \\
&= (|z|-|\zeta|)^2 + 2 |z|\,|\zeta|(1-\cos(|\arg \zeta - \arg z|)) \\
&\le \beta_j^2 + 2(a^j + \beta_j)^2 (1-\cos(2 C_0 a^{-j} \beta_j)) \\
&\le (1+4C_0^2(1+C_0^2)^2) \beta_j^2.
\end{split}\end{equation*}
Hence \eqref{eq:diameter-upper-bound} holds with $C_* = \sqrt{1+4C_0^2(1+C_0^2)^2}$.
\end{proof}

Set
$$
\Omega = \bigcup_{j \ge N} S_j.
$$
We call any open set $\Omega \subset \C$ constructed by this process a {\it sectorial sequence with scaling parameters $(b,a)$}. In the ensuing proofs, we write $A \lesssim B$, resp.\ $A \gtrsim B$, to indicate that two quantities $A(j)$ and $B(j)$ satisfy the inequality $A(j) \le C B(j)$, resp.\ $A(j) \ge (1/C)B(j)$ for some choice of $C>0$ and all $j \ge N$. We write $A \approx B$ if $A \lesssim B$ and $A\gtrsim B$.

\begin{proposition}
Let $\Omega$ be a sectorial sequence with scaling parameters $(b,a)$, and assume that $\beta_j = o(\eps^j)$ for every $\eps>0$. Then $\dim_L^\theta \Omega = 0$ for every $0<\theta<1$. In particular, $\dim_L \Omega = 0$.
\end{proposition}

For instance, we could choose $\beta_j = \exp(-c_1 e^{c_2j})$ for some $c_1,c_2>0$.

\begin{proof}
Assume that
\begin{equation}\label{eq:s-lower-bound}
\dim_L^\theta \Omega > s
\end{equation}
for some $0<\theta<1$ and $s>0$. We fix an integer $N_1 \gg N$ (to be determined later) and we restrict attention to $j \ge N_1$.

Fix $x \in S_j$ and let $R = \tfrac13 a^j$ and $r = C_* \beta_j$ where $C_*$ is the constant in \eqref{eq:diameter-upper-bound}. Then
$$
B(x,R) \cap \Omega = S_j
$$
in view of \eqref{eq:distance-lower-bound}. We claim that $r \le R^{1/\theta}$; the proof will be given soon. Since $r \ge \diam(S_j)$, the set $S_j$ can be contained in a single ball of radius $r$. Thus $1 \ge (1/C) (R/r)^s$ by \eqref{eq:s-lower-bound} whence
$$
a^j \lesssim \beta_j \lesssim b^j
$$
which contradicts the assumption that $b<a$.

It remains to show that $r \le R^{1/\theta}$. This inequality reads
\begin{equation}\label{eq:beta-to-a}
\beta_j \le C_*^{-1} 3^{-1/\theta} a^{j/\theta}.
\end{equation}
Since $\beta_j = o(\eps^j)$ for every $\eps>0$, there exists $N_1 \in N$ so that \eqref{eq:beta-to-a} holds true for all $j \ge N_1$. The proof is complete.
\end{proof}

In later applications we will be interested in sectorial sequences whose components $S_j$ are quasiround sets, with diameters tending to zero at a rate comparable to $b^j$. 

\begin{definition}
Let $\Omega = \cup_{j \ge N} S_j$ be a sectorial sequence open set with scaling parameters $(b,a)$ as defined above. We say that $\Omega$ is a {\bf special sectorial sequence} if $\beta_j \approx b^j$ and $\varphi_j \approx (b/a)^j$.
\end{definition}

\begin{proposition}
Let $\Omega$ be a special sectorial sequence with scaling parameters $(b,a)$. Set
$$
\theta_c := \frac{\log(1/a)}{\log(1/b)}.
$$
Then 
\begin{equation}\label{eq:lower-spectrum-formula}
\dim_L^\theta \Omega = \begin{cases} 2 \frac{\theta_c - \theta}{\theta_c(1-\theta)} & \mbox{if $0<\theta<\theta_c$,} \\ 0 & \mbox{if $\theta_c \le \theta < 1$.} \end{cases}
\end{equation}
In particular, $\dim_L \Omega = 0$.
\end{proposition}

\begin{proof}
First, we show that $\dim_L^\theta \Omega = 0$ if $\theta > \theta_c$. The proof is similar to the previous one; we only need to rework the proof of \eqref{eq:beta-to-a}. In the current situation, our assumption $\beta_j \approx b^j$ means that \eqref{eq:beta-to-a} is equivalent to $b^j \lesssim a^{j/\theta}$, and the latter condition holds true for $j \ge N_1$ and a sufficiently large choice of $N_1$ since $\theta > \theta_c$. By continuity of the lower spectrum in the parameter $\theta$, it also holds true that $\dim_L^{\theta_c} \Omega = 0$.

Now assume $0 < \theta < \theta_c$. Then $b > a^{1/\theta}$. We first show that $\dim_L^\theta \Omega \le 2 \tfrac{\theta_c-\theta}{\theta_c(1-\theta)}$. Fix $s$ with
\begin{equation}\label{eq:the-s}
0<s<\dim_L^\theta \Omega.
\end{equation}
For sufficiently large $j$, we may choose $x \in S_j$, $R = \tfrac13 a^j$ and $r = R^{1/\theta}$. Then $r = 3^{-1/\theta} a^{j/\theta} < 3^{-1/\theta} b^j$. The area of $S_j$ is
\begin{equation}\label{eq:area-of-Sj}
\int_{-\varphi_j}^{\varphi_j} \int_{a^j}^{a^j+\beta_j} r\,dr\, d\theta = \varphi_j ( 2 a^j \beta_j + \beta_j^2 ) \approx a^j \beta_j \varphi_j \approx b^{2j}.
\end{equation}
By area considerations, we may cover $B(x,R) \cap \Omega = S_j$ by $\approx (b^j/r)^2$ balls of radius $r$.  Due to \eqref{eq:the-s} we must have
$$
\left( \frac{b^j}{r} \right)^2 \gtrsim \left( \frac{R}{r} \right)^s
$$
or equivalently
$$
1 \gtrsim R^s r^{2-s} b^{-2j} \approx a^{j(s+(2-s)/\theta)} b^{-2j}.
$$
For this to hold true for all sufficiently large $j$ it must be the case that
\begin{equation}\label{eq:s-bound}
a^{s+(2-s)/\theta} \le b^2 \quad \Leftrightarrow \quad \left( s + \frac{2-s}{\theta} \right) \theta_c \ge 2 \quad \Leftrightarrow \quad s \le  2 \frac{\theta_c - \theta}{\theta_c(1-\theta)}.
\end{equation}
Letting $s \nearrow \dim_L^\theta \Omega$ we conclude that
$$
\dim_L^\theta \Omega \le 2 \frac{\theta_c - \theta}{\theta_c(1-\theta)}.
$$

Finally, we prove the opposite inequality $\dim_L^\theta \Omega \ge 2 \tfrac{\theta_c-\theta}{\theta_c(1-\theta)}$. To this end, let $x \in \Omega$, $0<r\le R^{1/\theta} < R$, and $s = 2 \tfrac{\theta_c-\theta}{\theta_c(1-\theta)}$. Fix $j$ so that $x \in S_j$ and recall that we are in the case $\theta < \theta_c$ or equivalently, $a< b^\theta$.

If $C_* \beta_j \le R \le \tfrac13 a^j$ then $B(x,R) \cap \Omega = S_j$ and
$$
r \le R^{1/\theta} \lesssim a^{j/\theta} < b^j \approx \beta_j.
$$
In view of the estimate \eqref{eq:area-of-Sj} for the area of $S_j$, we observe that we need $\approx (b^j/r)^2 \gtrsim (R/r)^s$ balls of radius $r$ to cover $S_j$.

If $R \le C_* \beta_j \lesssim b^j$ then $B(x,R) \cap \Omega \subset S_j$ and the area of $B(x,R) \cap S_j$ is comparable to $R^2$. In this case we need $\approx (R/r)^2$ balls of radius $r$ to cover $B(x,R) \cap S_j$. 

Finally if $R\ge\tfrac13 a^j$ then choose $\ell < j$ so that 
$$
\tfrac13 a^{\ell+1} \le R \le \tfrac13 a^\ell.
$$
Then $B(x,R) \cap \Omega \supset S_{\ell+1}$ and the area of $B(x,R) \cap \Omega$ is $\gtrsim a^{2\ell}$. In this case, we need $\gtrsim (a^\ell/r)^2\approx (R/r)^2$ balls of radius $r$ to cover $B(x,R) \cap \Omega$.

Since $s \le 2$ it follows in all cases that $B(x,R) \cap \Omega$ can be covered by $\gtrsim (R/r)^s$ balls of radius $r$.
We deduce that $\dim_L^\theta \Omega \ge s = 2 \tfrac{\theta_c-\theta}{\theta_c(1-\theta)}$ as claimed.
\end{proof}

\end{example}

\section{Quasiconformal distortion of lower dimension and the lower spectrum}\label{sec:qc-dist-lower}

This section is devoted to the proof of our main results on the quasiconformal distortion of lower dimension and lower spectrum. The conformal invariance of lower dimension is used in the proof of Theorem \ref{thm: QC L dimension distortion}. In section \ref{subsec:conf-invariance} we show that lower dimension is invariant under the inversion mapping
$\cI_a:\R^d \setminus \{a\} \to \R^d \setminus \{a\}$, $a \in \R^d$, defined as
\begin{equation}\label{eq:Ia}
\cI_a(x)
=
a+\frac{x-a}{|x-a|^2}
\end{equation}
under some auxiliary assumptions. This conclusion suffices for the required application in the proof of part (a) of Theorem \ref{thm: QC L dimension distortion}. Part (b) of that theorem then follows as a corollary.

\subsection{On the conformal invariance of lower dimension}\label{subsec:conf-invariance}

In this subsection we prove Proposition \ref{prop: dimL_inv} concerning the behavior of lower dimension under the conformal inversion $\cI_a$ defined in \eqref{eq:Ia}. We will later show that the conclusion of Proposition \ref{prop: dimL_inv} holds for arbitrary sets $E \subset \R^d \setminus \{a\}$. This conclusion will be established by using Proposition \ref{prop: dimL_inv} in conjunction with Theorem \ref{thm: QC L dimension distortion} (b) on global quasiconformal distortion of lower dimension. On the other hand, Proposition \ref{prop: dimL_inv} itself will play a key role in the proof of Theorem \ref{thm: QC L dimension distortion}(b).

\begin{proposition}
\label{prop: dimL_inv}
Fix $a \in \R^d$, and let $E\subset \R^d \setminus \{a\}$ satisfy either
\begin{itemize}
\item[(a)] $E$ is closed and unbounded, or
\item[(b)] $E$ is bounded and $a \in \partial E$.
\end{itemize}
Then $\dim_L \cI_a(E)=\dim_L E$.
\end{proposition}

\begin{proof}
Since $\cI_a$ is an involution and $E$ is closed and unbounded if and only if $\cI(E)$ is bounded with $a \in \partial \cI_a(E)$, it suffices to consider only one of the two cases. We assume for the remainder of the proof that $E$ is closed and unbounded.

Since the lower dimension is invariant under similarities, it suffices to assume $a=0$. Set $\cI=\cI_0$ and $F=\cI(E)\cup \{0\}$. Note that by invariance of the lower dimension under closure, we have 
$$
\dim_L F=\dim_L \cI(E),
$$ 
which allows us to focus on the dimension of the closed set $F$.
Since $E$ is closed, we have
$$
\delta:=\dist(0,E)>0,
$$
and there exists $x_0\in E$ with $|x_0|=\delta$.
Hence, $F\subset B(0,\delta^{-1})$ and
$$
\delta^{-1}\leq\diam F\leq 2\delta^{-1},
$$ due to $\cI(x_0), 0 \in F$ and $|\cI(x_0)-0|=\delta^{-1}$.

We  prove that
$$
\dim_L F\geq\dim_L E,
$$ and the opposite inequality can be proved similarly. If $\dim_L E=0$ the inequality is trivially true, so suppose $\dim_L E>0$ and let $s\in (0,\dim_LE)$.
By definition of the lower dimension, there is a constant $C>0$ such that
\begin{equation}\label{eq: lower-estimate-E-inversion}
N_r(B(x,R)\cap E)
\geq
(1/C) \left(\frac{R}{r}\right)^s
\end{equation}
for every $x\in E$ and $0<r<R<\infty$. Fix $L>1$ so that $L^s/C\geq 2$.

Let $y\in F$ and $0<r'<R'<\diam F$. The goal is to find a uniformly small ball in $B(y,R')$ that is far enough from $0$, so that applying $\cI$ on the non-empty intersection of said ball with $F$ yields a ball intersecting $E$ in the source, where we can apply \eqref{eq: lower-estimate-E-inversion}. 

\begin{claim}\label{claim: INV}
    There is $c>0$, independent of $y, r', R'$, 
and a point $w\in F$ such that
\begin{equation}\label{eq: ball-ByR'-far-0}
    cR'\leq \frac{|w|}{4} \quad \text{and} \quad B(w,cR')\subset B(y, R').
\end{equation}
\end{claim}

\begin{proof}[Proof of Claim \ref{claim: INV}]
    Let $c:=(200L)^{-1}>0$. If $|y|\geq R'/10$, then $w=y$ satisfies \eqref{eq: ball-ByR'-far-0}, so suppose $|y|<R'/10$. Applying \eqref{eq: lower-estimate-E-inversion} for $x=x_0$, $r=20/R'$, and $R=20L/R'$ yields
$$
N_{20/R'}(B(x_0,20L/R')\cap E)\geq L^s/C\geq 2,
$$ the last inequality due to the choice of $L$. Hence, $B(x_0,20L/R')\cap E$ cannot be contained in $B(x_0,10/R')$ and, therefore, there is some $z\in B(x_0,20L/R')\cap E$ with
$$
\frac{10}{R'}<|z-x_0|\leq \frac {20L}{R'}.
$$ Since $|x_0|=\delta$ and $R'<\diam F\leq 2\delta ^{-1}$, also noting that $10/R'-\delta>8/R'$, the above yields
$$
\frac{R'}{2(10L+1)}\leq \frac 1{|z|}\leq \frac{R'}8.
$$ Hence, $w=\cI (z)$ satisfies \eqref{eq: ball-ByR'-far-0} due to $|w|=|z|^{-1}$, the choice of $c$, the assumption $|y|<R'/10$, and the above inequality.
\end{proof}

Let $w\in F$ as in \eqref{eq: ball-ByR'-far-0} and let $z=\cI^{-1}(w)$. Set
\begin{equation}\label{eq: INV pair r-R}
    r:= 2r'|z|^2, \qquad R:= \frac{cR'|z|^2}{4}.
\end{equation} If $u\in B(z, R)$ then by \eqref{eq: ball-ByR'-far-0} and $|z|=|w|^{-1}$ we have
$$
|u|\geq |z|-\frac{cR'|z|^2}{4}\geq |z|-\frac{|z|}{16}=\frac{15|z|}{16}.
$$ 
The above and an elementary property of $\cI$ yield
$$
|\cI(u)-w|=\frac{|u-z|}{|u|\,|z|}\leq \frac{cR'|z|^2}{4|u|\,|z|}\leq \frac{4cR'}{15}< cR'.
$$ Since $u$ is arbitrary, the above and \eqref{eq: ball-ByR'-far-0} imply
$$
\cI(B(z,R))\subset B(w, cR')\subset B(y,R').
$$ Therefore, if $\{V_i\}_{i\in I}$ is a collection of sets of diameter at most $r'$ that covers $\cI(B(z,R)\cap E)$, the collection $\{\cI(V_i)\}_{i\in I}$ covers $B(z,R)\cap E$. Suppose that $v, v'\in V_i$ have images under $\cI$ contained in $B(z,R)$. Due to $\cI(v)\in B(z,R)$, \eqref{eq: ball-ByR'-far-0} and $|z|=|w|^{-1}$ we have
$$
|\cI(v)|\leq R+|z|=\frac{cR'|z|^2}{4}+|z|\leq \frac{|z|}{16}+|z|=\frac{17|z|}{16},
$$ and similarly for $\cI(v')$. This implies that
$$
|\cI(v)-\cI(v')|=\frac{|v-v'|}{|v|\,|v'|}\leq 2r' |z|^2.
$$ Since $v,v'$ are arbitrary, the above yields
$$
\diam (\cI(V_i)\cap B(z,R))\leq r,
$$ for all $i\in I$. Since $\{V_i\}_{i\in I}$ is an arbitrary collection of sets of diameter at most $r'$ that covers $\cI(B(z,R)\cap E)\subset B(y,R')\cap F$, we have
\begin{equation}\label{eq: INV_cover_numbers}
    N_{r'}(B(y,R')\cap F)\geq N_r(B(z,R)\cap E).
\end{equation}
If $r\geq R$, then by \eqref{eq: INV pair r-R} we have
$$
\left(\frac{R'}{r'}\right)^s\leq \left( \frac8c\right)^s,
$$ which implies
\begin{equation}\label{eq: INV N-image-C1}
    \left( \frac c8\right)^s\left(\frac{R'}{r'}\right)^s\leq 1\leq N_{r'}(B(y,R')\cap F).
\end{equation}
Thus, suppose $r<R$. This allows us to use \eqref{eq: lower-estimate-E-inversion}, which along with \eqref{eq: INV_cover_numbers} imply
\begin{equation}\label{eq: INV N-image-C2}
N_{r'}(B(y,R')\cap F)\geq (1/C) \left( \frac{R}{r} \right)^s=\frac{c^s}{8^s C}\left(\frac{R'}{r'}\right)^s.
\end{equation} Therefore, by \eqref{eq: INV N-image-C1} and \eqref{eq: INV N-image-C2} we have
$$
N_{r'}(B(y,R')\cap F)\geq (1/\tilde{C}) \left(\frac{R'}{r'}\right)^s,
$$ for all $y\in F$ and all $r'<R'<\diam F$, where $\tilde{C}=\max \{\frac{8^s C}{c^s}, \frac{8^s}{c^s}\}$. This implies that
$$
\dim_L F\geq s,
$$ 
and letting $s\to \dim_L E$ completes the proof.

For the reverse inequality, without loss of generality assume that
$\dim_L F>0$, and fix
$s\in(0,\dim_L F)$. There is a constant $C>0$ such that
\begin{equation}\label{eq: lower-estimate-F-inversion}
N_{r'}(B(y,R')\cap F)
\geq
(1/C) \left(\frac{R'}{r'}\right)^s
\end{equation}
for every $y\in F$ and every
$0<r'<R'<\diam F$. Fix $L>1$ so that $L^s/C\geq2$.

Similarly to the previous direction, we make the following claim.

\begin{claim}\label{claim: INV2}
For every $x\in E$ and every $R>0$, there is a point $z\in E$
such that
\begin{equation}\label{eq: ball-BxR-far-0}
\frac{R}{200L}\leq\frac{|z|}{4}
\qquad\text{and}\qquad
B(z,R/(200L))\subset B(x,R).
\end{equation}
\end{claim}

\begin{proof}[Proof of Claim \ref{claim: INV2}]
The proof is almost identical to that of Claim \ref{claim: INV}.
\end{proof}

Fix $x\in E$ and $0<r<R<\infty$, and choose $z\in E$ as in Claim
\ref{claim: INV2}. Let $y:=\cI(z)$
and define
\begin{equation}\label{eq: INV-pair-r-R-reverse}
r':=2r|y|^2,
\qquad
R':=\frac{R|y|^2}{1000L}.
\end{equation}
Note that unlike the previous direction, $R'<\infty$ is not enough and the inequality $R'<\diam F$ in needed in order to apply \eqref{eq: lower-estimate-F-inversion}. By $|y|=|z|^{-1}$ and \eqref{eq: ball-BxR-far-0} we have
$$
\frac{R|y|}{200L}\leq\frac14,
$$
which along with $0,y\in F$ and \eqref{eq: INV-pair-r-R-reverse} implies
$$
R'\leq\frac{|y|}{20}<\diam F.
$$ The rest of the proof follows similarly to the proof of the opposite inequality.
\end{proof}

\subsection{Proof of lower spectrum distortion}

Recall that
$$
\theta(t) = \frac{1}{1+t}.
$$
For $t>0$, $K \ge 1$, and $d \ge 2$ we also set
$$
\ut = \ut(d,K) := \left(1 - \frac{d}{p_O^{\RH}(d,K)} \right)^{-1} t
$$
and
$$
\ot = \ot(d,K) := \left(1 - \frac{d}{p_I^{\RH}(d,K)} \right) t.
$$
We shall prove a slightly stronger estimate that involves the reverse H\"older inequality exponent that depends on the inner dilatation, namely the exponent $p_I^{\RH}(d,K)$. In particular, due to $p_O^{\RH}(d,K^{d-1})\leq p_I^{\RH}(d,K)$ and the fact that $K_I(f^{-1})=K_O(f)$, the lower estimate in the following theorem implies the one in Theorem~\ref{thm: QC L spectrum distortion}.

\begin{theorem}\label{thm:KI QC L spectrum distortion}
Let $f:\Omega\rightarrow\Omega'$ be a $K$-quasiconformal mapping between domains in $\R^d$, $d \ge 2$. Suppose $t\in (0,\infty)$ and $E\subset \Omega$ is compact. Then
\begin{equation}\label{eq:KI dim-l-spec}
\frac{(p_I^{\RH}(d,K)-d)\,\dim_{L,\reg}^{\theta(\ot/K^{d-1})} E}
{p_I^{\RH}(d,K)-\dim_{L,\reg}^{\theta(\ot/K^{d-1})} E}
\le \dim_{L,\reg}^{\theta(t)} f(E)
\le
\frac{p_O^{\RH}(d,K)\,\dim_{L,\reg}^{\theta(K \ut)} E}
{p_O^{\RH}(d,K)-d+\dim_{L,\reg}^{\theta(K \ut)} E}.
\end{equation}
\end{theorem}

\begin{proof}
It suffices to prove the right hand side inequality in \eqref{eq:KI dim-l-spec}, as the left hand side follows by applying the upper bound to the $K^{d-1}$-quasiconformal mapping $f^{-1}$.
    
Let $t>0$, $\alpha>\alpha_E$, $p<p^{\RH}_O(d,K)$, and set $\beta=\frac{p \alpha}{p-d+\alpha}$ and $\widetilde{t_p}=(1-d/p)t$. By definition of $\dim_{L,\text{reg}}^{\theta(Kt)} E$ we have that for every $n\in \N$ there are $0<2^{-m_n} R_n\leq R_n^{1/\theta(Kt)}<R_n<\diam E$ and $x_n\in E$ so that
    \begin{equation}\label{eq: def of dimL E}
        N(B(x_n,R_n)\cap E, m_n)< \frac{1}{n} 2^{m_n \alpha}.
    \end{equation} In the above we used the definition of lower spectrum with dyadic cubes. Set $Q_n:=Q(x_n,R_n)$. After potentially passing to a subsequence of $(m_n)_{n\in\N}$, we may assume that $m_n\rightarrow\infty$ and that $Q(x_n,2R_n)\subset \Omega$ for all $n\in\N$. Moreover, we further reduce the sequence $(m_n)_{n\in \N}$ only to those terms greater or equal to the smallest positive $m_{n_0}\in \N$ so that
    \begin{equation}\label{eq: Qi away partial Omega}
        \diam Q_{n_0} \leq \frac{\dist (B(E,\sqrt{d} \,R_{n_0}),\partial \Omega)}{10}
    \end{equation}
    and 
    \begin{equation}\label{eq: f(Qi) away partial Omega}
        (\diam Q_{n_0})^{1-d/p} \leq \frac{\dist(B(f(E),L_{n_0}),\partial \Omega')}{10 C_{MS} (\int_{B(E,2^{-1}d(E,\partial \Omega)) }|Df|^p)^{1/p}},
    \end{equation}
    where $L_n= C_{MS} (\diam Q_n)^{1-d/p}(\int_{B(E,2^{-1}d(E,\partial \Omega)) }|Df|^p)^{1/p}$ and $C_{MS}$ is the Morrey-Sobolev constant (see \eqref{eq:morreysobolev}). Note that such a choice for $n_0$ in \eqref{eq: Qi away partial Omega} is possible due to the left-hand-side converging to $0$ as $n_0$ increases, while the right-hand-side converges to a positive uniform constant, and similarly for the choice in \eqref{eq: f(Qi) away partial Omega}. The inequalities \eqref{eq: Qi away partial Omega}, \eqref{eq: f(Qi) away partial Omega} essentially ensure that we may apply the reverse H\"older inequality for $p$ on the cubes $Q_n$ in following arguments. Furthermore, by bi-Lipschitz stability of the lower spectrum, we may assume henceforth that $\diam E, \diam f(E)<1$.
    
    Let $\{Q_i^{m_n}\}_{i\in I_{m_n}}$ be a minimal collection of $m_n$-level dyadic sub-cubes of $Q(x_n,R_n)$ covering $B(x,R_n)\cap E$. If for all $n\in \N$ we find $y_n\in f(E)$, $0<r_n'\lesssim (R_n')^{1/\theta(t)}<R_n'<\diam f(E)$ so that
    \begin{equation}\label{eq: needed covering dimL f(E)}
        N_{r_n'}(B(y_n,R_n')\cap f(E))<a_n \left(\frac{R_n'}{r_n'}\right)^\beta,
    \end{equation} for some positive sequence $a_n\rightarrow 0$, 
    then we have $\dim_{L,\text{reg}}^{\theta(t)} f(E)\leq \beta$ by definition of lower spectrum. Taking $\alpha\rightarrow \alpha_E$ and $p\rightarrow p^{\RH}_O(d,K)$ after that finishes the proof.

    Hence, we focus on finding such sequences of points and scales in the target. For $n\in \N$, let $y_n=f(x_n)$ and
    $$
    R_n'=\frac{\inf \{|y_n-f(x)|: \,x\in \partial B(x_n,R_n)\}}{2}.
    $$ This choice implies that $B(y_n,R_n')\subset f(B(x_n,R_n))$, and it is evident that covering 
    $$
    f(B(x_n,R_n))\cap f(E)
    $$ by images under $f$ of cubes $Q_i^{m_n}$ and/or their sub-cubes  results in a covering of $B(y_n,R_n')\cap f(E)$. Also, set
    $$
    r_n'= 2^{-\tfrac{m_n \alpha}{\beta}} R_n'.
    $$ Note that by $r_n=2^{-m_n}R_n\leq R_n^{1/\theta(Kt)}$ and $R_n\lesssim(R_n')^{1/K}$, the latter due to H\"older continuity of $f^{-1}$, we have that 
    \begin{equation}\label{eq: towards rn' admissible}
        r_n'\leq R_n^{(1/\theta(Kt)-1)\beta^{-1} \alpha}R_n' \lesssim (R_n')^{K^{-1}(1/\theta(Kt)-1)\beta^{-1} \alpha+1}.
    \end{equation} However, by $\theta(t)=1/(1+t)$ and choice of $\beta$, we have that
    $$
    K^{-1}(1/\theta(Kt)-1)\beta^{-1} \alpha+1=\beta^{-1} \alpha \,t+1=\frac{p-d+\alpha}{p}\,t+1\geq 1/\theta(\widetilde{t_p}),
    $$ where $\widetilde{t_p}=(1-d/p)t$. By the assumption that $\diam f(E)<1$ and by \eqref{eq: towards rn' admissible}, the above implies
    $$
    r_n'\lesssim (R_n')^{1/\theta(\widetilde{t_p})},
    $$ as needed. Thus, the small scales $r_n'$ are admissible for the $\theta(\widetilde{t_p})$-lower spectrum of $f(E)$.

    Since $p<p^{\RH}_O(d,K)$, by \eqref{eq: Qi away partial Omega}, \eqref{eq: f(Qi) away partial Omega} and an application of Morrey-Sobolev inequality on $Q_n$, we conclude that the reverse H\"older inequality is satisfied on all $Q_n$, i.e.
    $$
    \frac{1}{|Q(x_n,R_n)|^{1/p}} \left( \int_{Q(x_n,R_n)} |Df|^p \right)^{1/p} \leq \frac{C_{\RH}}{|Q(x_n,R_n)|^{1/d}}\left( \int_{Q(x_n,R_n)} |Df|^d \right)^{1/d} \, .
    $$ 
    This implies that
    $$
    \frac{1}{R_n^d}  \int_{Q(x_n,R_n)} |Df|^p  \lesssim \frac{1}{R_n^p}\left( \int_{Q(x_n,R_n)} J_f \right)^{p/d} \, ,
    $$which in turn leads to
    \begin{equation}\label{eq: working RHI estimate}
        R_n^{p-d} \int_{Q(x_n,R_n)} |Df|^p\lesssim |f(Q(x_n,R_n)|^{p/d}.
    \end{equation} However, by Proposition~\ref{prop:egg-yolk}, $f$ is quasisymmetric on $B(E, 2\sqrt d R_{n_0})$, which by choice of $R_n'$ yields that 
    $|f(Q(x_n,R_n)|\lesssim (R_n')^d$.
    
    Combining this with \eqref{eq: working RHI estimate} implies
    \begin{equation}\label{eq: RHI big scales relation}
        R_n^{p-d}(R_n')^{-p}\lesssim \left( \int_{Q(x_n,R_n)} |Df|^p \right)^{-1}.
    \end{equation}

    Fix $m=m_n$ and split cubes $Q_i^{m}$ of this level between two different collections. We say that $Q_i^{m}$ is \textit{$r_m'$-major} if $\diam f(Q_i^{m})> r_m'$, and we say that $Q_i^{m}$ is \textit{$r_m'$-minor} otherwise. If all cubes $Q_i^{m}$ are $r_m'$-minor, we have by \eqref{eq: def of dimL E} that
    \begin{equation}\label{eq: all minor}
    N_{r_n'}(B(y_n,R_n')\cap f(E))<\frac{1}{n} 2^{m_n\alpha}=\frac{1}{n}\left(\frac{R_{m_n}'}{r_{m_n}'}\right)^\beta.
    \end{equation} Suppose now that there are $r_m'$-major cubes $Q_i^{m}$. Note that after sub-dividing every $m$-level $r_m'$-major cube, there might still be $r_m'$-major cubes among the children, i.e. among $(m+1)$-level cubes, and so on. We denote by $\{Q_i^{k}\}_{i\in J_{m,k}}$ the collection of all $r_m'$-major cubes of level $k\geq m$, and by $M(k)$ their number, i.e. $M(k)=\# J_{k,m}$.
    The Morrey--Sobolev inequality on a major $Q^k_i$ takes the form
    $$
    \diam f(Q_i^k) \leq C_{MS} (\diam (Q_i^k))^{1-d/p} \left( \int_{Q_i^k} |Df|^p \right)^{1/p}.
    $$
    Since $Q_i^k$ is $r_m'$-major, this implies that
    $$
    2^{-\frac{m\alpha p}{\beta}}(R_n')^p \leq C_{MS}^p \, 2^{-k (p-d)} R_n^{p-d} \int_{Q_i^k}|Df|^p
    $$
    from which, by summing over all $r_m'$-major cubes of level $k$, it follows that
    \begin{equation*}
    2^{-\frac{m\alpha p}{\beta}}(R_n')^p M(k)
    \leq C_{MS}^p \, 2^{-k(p-d)} R_n^{p-d} \int_{\bigcup\limits_{i\in J_{k,m}} Q_i^k}|Df|^p .
    \end{equation*}
    Summing over all $k \geq m$ yields
    $$
    2^{-\frac{m\alpha p}{\beta}}(R_n')^p \sum_{j=m}^{\infty} M(j) \leq \tilde{C}_{Ms} \, 2^{-m(p-d)} R_n^{p-d} \int_{\bigcup\limits_{k\geq m}\,\bigcup\limits_{i} Q_i^k} |Df|^p,
    $$ where $ \tilde{C}_{Ms}$ is   $C_{MS}^p$ multiplied by a uniform constant that is due to the elementary property of geometric series starting at a term $m$.
    Hence, by setting
    $$
    D_n=\frac{\int_{\bigcup\limits_{k\geq m}\,\bigcup\limits_{i} Q_i^k}|Df|^p}{\int_{Q(x_n,R_n)} |Df|^p},
    $$
    and using \eqref{eq: RHI big scales relation}, we have
    \begin{equation}\label{eq: total major cubes step before}
    \sum_{k=m}^{\infty} M(k) \lesssim \tilde{C}_{MS} C_{RH} \, 2^{\frac{m\alpha p}{\beta}-m(p-d)} \, D_n \,.
    \end{equation}
    However, by choice of $\beta$ we get
    \begin{equation}\label{eq: total major cubes}
    \sum_{k=m}^{\infty} M(k) \leq C \, 2^{m\alpha} \, D_n \,,
    \end{equation}
    where $C>0$ is the uniform constant resulting from combining $C(\lesssim)$ in \eqref{eq: total major cubes step before} and $\tilde{C}_{MS} C_{RH}$. If we subdivide all $r_m'$-major cubes of all levels, we are left with at most $2^d\sum_{k=m}^{\infty} M(k)+\#I_m$ $r_m'$-minor cubes. Therefore, we have by \eqref{eq: total major cubes} that
    \begin{equation}\label{eq: last covering rel proof}
        N_{r_n'}(B(y_n,R_n')\cap f(E))\leq C \, 2^{m\alpha} \, D_n+ \frac{1}{n} 2^{m\alpha}.
    \end{equation}

    Set $E_n:=\bigcup\limits_{k\geq m}\,\bigcup\limits_{i} Q_i^k=\bigcup\limits_{k\geq m_n}\,\bigcup\limits_{i} Q_i^k$, then by \eqref{eq: def of dimL E}, the quotient of $d$-Hausdorff measures
    \begin{equation}\label{eq: quotient volumes 0}
        \frac{|E_n|}{|Q(x_n,R_n)|}\lesssim \, \frac{\tfrac{1}{n}2^{m_n \alpha}(2^{-m_n}R_n)^d}{R_n^d}=\frac{1}{n2^{m_n(d-\alpha)}}\rightarrow 0
    \end{equation} as $n\rightarrow\infty$. We show that this convergence is enough to guarantee that $D_n\rightarrow 0$. 
    Fix $r\in (0,1)$ so that $p(1+r)<p^{\RH}_O(d,K)$. By H\"older inequality and the inclusion $E_n\subset Q_n$ we have
    $$
    \int_{E_n} |Df|^p\leq \left( \int_{E_n} |Df|^{p(1+r)} \right)^{\tfrac{1}{1+r}} |E_n|^{\tfrac{r}{1+r}}\leq \left( \int_{Q_n} |Df|^{p(1+r)} \right)^{\tfrac{1}{1+r}} |E_n|^{\tfrac{r}{1+r}},
    $$
    which after an application of the reverse H\"older inequality and another application of H\"older's inequality yields
    \begin{align*}
        \int_{E_n} |Df|^p
        &\lesssim
        |Q_n|^{\tfrac{1}{1+r}}
        |Q_n|^{-\tfrac{p}{d}}
        \left( \int_{Q_n} |Df|^d \right)^{p/d}
        |E_n|^{\tfrac{r}{1+r}}
        \\
        &\le
        |Q_n|^{\tfrac{1}{1+r}-\tfrac{p}{d}}
        |E_n|^{\tfrac{r}{1+r}}
        \left( \int_{Q_n} |Df|^p \right)
        |Q_n|^{\tfrac{p-d}{d}} .
    \end{align*}
    Hence, the above implies
    $$
    D_n\lesssim \left( \frac{|E_n|}{|Q_n|} \right)^{\tfrac{r}{1+r}},
    $$ which by \eqref{eq: quotient volumes 0} results in $D_n\rightarrow0$ as $n\rightarrow \infty$. 
    Therefore, the sequence
    $$
    a_n:= 100(C \,D_n+ \frac{1}{n})
    $$ is enough to satisfy \eqref{eq: needed covering dimL f(E)}, potentially after passing to a subsequence. This completes the proof of the right-hand side inequality, after taking $p\to p^{\RH}_O(d,K)$, $\alpha\to \alpha_E$ and using the continuity of the regularized lower spectrum.  The left-hand side inequality follows by applying the right-hand side bound on $f^{-1}$ and $K_I(f^{-1})=K_O(f)$. 
\end{proof}

\begin{remark}
The corresponding result for Assouad spectrum was proved in \cite{OurQCspec}. A slightly different approach was taken there, namely, the diameter conditions necessary for application of the reverse H\"older inequality were obtained by partitioning the compact set $E$ via a Whitney decomposition of the domain and using finite stability of the Assouad spectrum. Since the lower spectrum (and lower dimension) are not finitely stable, a new approach is needed here. In lieu of the Whitney decomposition argument we rely on the fact that the necessary conditions \eqref{eq: Qi away partial Omega} and \eqref{eq: f(Qi) away partial Omega} hold for cubes $Q_n$ with $n$ sufficiently large.
\end{remark}

\begin{remark}\label{rem: spectrum exp diff}
    Compared to the regularized Assouad spectrum distortion theorem in \cite{OurQCspec}, the admissibility spectrum exponents are different. This is due to the flexibility in the Assouad spectrum case of going down enough levels to achieve admissibility in the considered pairs of scales, while this is not an option in the regularized lower spectrum case. In particular, the pairs of scales are already fixed by the choice of the sequence $1/n$ corresponding to $\alpha>\alpha_E$ and there is no control over the ``gaps'' between consecutive levels of pairs $(r_n, R_n)$.
\end{remark}

\begin{remark}
Fraser and Yu \cite{fy:assouad-spectrum} (see also Lemma 3.4.13 in \cite{FraserBook}) studied the distortion of the lower spectrum by bi-H\"older homeomorphisms. Since quasiconformal maps are locally bi-H\"older, it is instructive to consider the relationship between Theorem \ref{thm: QC L spectrum distortion} and the results of \cite{fy:assouad-spectrum}. Our discussion in this remark parallels our previous discussion of the Assouad spectrum in\cite[pages 302--303]{OurQCspec}.

Recall that a homeomorphism $f:X\to Y$ between metric spaces is said to be {\it $(\alpha,\beta)$-bi-H\"older}, for $0<\alpha\le 1\le\beta <\infty$, if there exists a constant $C>0$ so that
$$
C^{-1}d(x,x')^\beta \le d(f(x),f(x')) \le C d(x,x')^\alpha \qquad \forall\,x,x' \in X.
$$
If $\alpha < \beta$, it is clear that $X$ must be a bounded space in order for a $(\alpha,\beta)$-bi-H\"older homeomorphism from $X$ to exist. According to \cite[Lemma 3.4.13]{FraserBook},
if $X \subset \R^n$ is bounded, $f:X \to \R^n$ is a $(\alpha,\beta)$-bi-H\"older homeomorphism, and $0<\theta<\alpha/\beta$, then
\begin{equation}\label{eq:biHolder}
\frac{1-\beta\theta/\alpha}{\beta(1-\theta)} \dim_{L,reg}^{\beta\theta/\alpha}(X) \le \dim_{L,reg}^\theta(f(X)) \le \frac{1-\alpha\theta/\beta}{\alpha(1-\theta)} \dim_{L,reg}^{\alpha\theta/\beta}(X).
\end{equation}
Every $K$-quasiconformal map in $\C$ is locally $(\tfrac1K,K)$-bi-H\"older continuous, see e.g.\ \cite[Corollary 3.10.3]{aim:2qc}. It follows from \eqref{eq:biHolder} that if $f:\C\to\C$ is $K$-quasiconformal, $E\subset\C$ is bounded, and $0<\theta<1/K^2$, then
\begin{equation}\label{eq:biHolder2}
\dim_{L,reg}^\theta(f(E)) \le K\, \frac{1-\theta/K^2}{1-\theta} \dim_{L,reg}^{\theta/K^2}(E).
\end{equation}
On the other hand, \eqref{eq: dim-spec-dimension-2} implies that if $t>0$ then
\begin{equation}\label{eq:ours}
\dim_{L,reg}^{\theta(t)}(f(E)) \le K \frac{\dim_{L,reg}^{\theta(K^2 t)}(E)}{1+\tfrac{K-1}{2}\dim_{L,reg}^{\theta(K^2 t)}(E)} \,.
\end{equation}
The following observation confirms that the estimate in \eqref{eq:ours} provides a strict improvement over that in \eqref{eq:biHolder2}, for all choices of $t$ so that the necessary spectrum parameters are defined.

\begin{proposition}\label{prop:final}
Let $d:= \dim_{L,reg}^{\theta(K^2t)}(E)$. If $\theta(t)<1/K^2$, then
$$
\frac{d}{1+\tfrac{K-1}{2}d} \le \frac{1-\theta(t)/K^2}{1-\theta(t)} \dim_{L,reg}^{\theta(t)/K^2}(E).
$$
\end{proposition}

\begin{proof}[Proof of Proposition \ref{prop:final}]
Let $\theta_1 = \theta(t)/K^2$ and $\theta_2 = \theta(K^2t)$. Then $d = \dim_{A,reg}^{\theta_2}(E)$ and the conclusion reads
\begin{equation}\label{eq:prop-conclusion}
\frac{d}{1+\tfrac{K-1}{2}d} \le \frac{1-\theta_1}{1-\theta(t)} \dim_{L,reg}^{\theta_1}(E).
\end{equation}
Since $\theta_1 \le \theta_2$, monotonicity implies that
$$
\dim_{L,reg}^{\theta_1}(E) \ge \dim_{L,reg}^{\theta_2}(E)
$$
and so \eqref{eq:prop-conclusion} is implied by
\begin{equation}\label{eq:prop-conclusion2}
\frac{d}{1+\tfrac{K-1}{2}d} \le \frac{1-\theta_1}{1-\theta(t)}d.
\end{equation}
Since $1-\theta_1=\tfrac{K^2(1+t)-1}{K^2(1+t)}$ and $1-\theta(t)=\tfrac{t}{1+t}$, \eqref{eq:prop-conclusion2} reads 
$$
\frac1{1+\tfrac{K-1}{2}d} 
\le 
1 + \frac{K^2-1}{K^2 t},
$$
which is satisfied for all $t>0$ and $K\ge 1$.
\end{proof}
\end{remark}

\subsection{Proof of lower dimension distortion}

Similarly to the previous subsection, we shall prove the following stronger theorem involving the inner dilatation exponent, which implies Theorem~\ref{thm: QC L dimension distortion}. 

\begin{theorem}\label{thm:KI QC L dimension distortion}
Let $f:\Omega\rightarrow\Omega'$ be a $K$-quasiconformal map between domains in $\R^d$, $d \ge 2$.

\medskip

\noindent (a) For any compact set $E\subset \Omega$, it holds true that
\begin{equation}\label{eq:KI dimL}
 \frac{(p_I^{\RH}(d,K)-d)\,\dim_{L} E}{p_I^{\RH}(d,K)-\dim_{L} E}\leq \dim_{L} f(E) \le \frac{p_O^{\RH}(d,K)\,\dim_{L} E}{p_O^{\RH}(d,K)\, - d + \dim_{L} E}
\end{equation}
(b) If $\Omega = \Omega' = \R^d$, then the estimate in \eqref{eq:KI dimL} holds true for arbitrary sets $E \subset \R^d$.
\end{theorem}

\begin{proof}
(a) Similarly to Theorem \ref{thm: QC L spectrum distortion}, it is enough to prove the upper bound in \eqref{eq:KI dimL}. The proof follows verbatim from the proof of Theorem~\ref{thm: QC L spectrum distortion} by replacing the two spectrum variables $\theta(Kt)$ and $\theta((1-d/p)t)$ by the constant $1$, with the admissibility of the $r_n'$ scales following directly in this case.

(b)We now prove that the compactness assumption can be removed in case the $K$-quasiconformal mapping $f$ is globally defined.

Let $E$ be an arbitrary subset of $\R^d$. If $E = \R^d$ the conclusion is trivial, so assume $E \ne \R^d$. Without loss of generality we may assume that $E$ is closed. Choose $a \not \in E$ and consider the map
$$
g = \cI_{f(a)} \circ f \circ \cI_a.
$$
Then $g:\R^d \setminus \{a\} \to \R^d \setminus \{f(a)\}$ is a $K$-quasiconformal mapping. Isolated boundary points are removable singularities of quasiconformal mappings, see for instance \cite[Theorem 17.3]{Vai:book}. We may thus extend $g$ to a $K$-quasiconformal map (also denoted $g$) from $\R^d$ to $\R^d$. Now $\cI_a(E)$ is a precompact subset of $\R^d$ and we may consider its compact closure $\overline{\cI_a(E)} = \cI_a(E) \cup \{a\}$ which has the same lower dimension. The proof is concluded by an application of Theorem \ref{thm: QC L spectrum distortion}(a) to the set $\overline{\cI_{a}(E)}$. 
\end{proof}

\begin{remark}
Passing from compact sets to arbitrary sets in the distortion estimates under quasiconformal and conformal mappings is simpler in the Assouad dimension case, since Assouad dimension is finitely stable. The more involved approach adopted in this paper is required to overcome the lack of finite stability for lower dimension.
\end{remark}

\begin{corollary}
Globally defined conformal mappings of $\R^d$ preserve the lower dimension of arbitrary subsets.
\end{corollary}

\begin{proof}
In all dimensions $d \ge 2$, every globally defined conformal mapping is M\"obius. It suffices to show that inversions preserve the lower dimensions of arbitrary subsets. The invariance of lower dimension under the inversion $\cI_a$ was proved in Proposition \ref{prop: dimL_inv} for sets $E \subset \R^d \setminus \{a\}$ which are either closed and unbounded, or bounded with $a$ in the closure. These conclusions were then used to prove Theorem 
\ref{thm: QC L dimension distortion}(b). The proof is completed by an appeal to Theorem \ref{thm: QC L dimension distortion}(a), noting that conformal mappings are $1$-quasiconformal and that $p^{\RH}(d,1) = \infty$ (see for instance \cite[Proposition 2 c)]{Iwan_Lp_int_in_PDEs}).
\end{proof}

\begin{remark}\label{rem:why-reciprocals}
We discuss in more detail the role of the symmetric estimates in \eqref{eq: dim-l-1} and \eqref{eq: dim-l-2} in the theory of quasiconformal dimension distortion. In addition the visual appeal of these symmetric upper and lower bounds, recent developments suggest that reciprocal dimension, and more precisely the `relative reciprocal dimension'
$$
\frac1{\dim E} - \frac1d
$$
of a set $E \subset \R^d$ is a natural quantity to consider. 

First, the recent work of Fuhrer--Ransford--Younsi \cite{fry:holomorphic-motions}, Menssen--Younsi \cite{my:holomorphic-motions} and Ransford--Younsi \cite{ry:holomorphic-motions} shows that reciprocal dimension is often an inf-harmonic quantity under planar holomorphic motions. More specifically, the function $\lambda \mapsto 1/\dim(E_\lambda) - 1/2$ is inf-harmonic when $E_\lambda$ varies holomorphically in $\C$, when $\dim$ is either upper box-counting dimension, packing dimension, or quasi-Assouad dimension. However, inf-harmonicity fails if $\dim$ is either Hausdorff dimension or lower box-counting dimension. The validity of the sharp dimension distortion estimates then follows as an immediate consequence of the inf-harmonicity.

Second, in \cite{FraTyson_intermed_dim}, quasiconformal distortion of intermediate dimensions (another notion of interpolation which interpolates between Hausdorff and upper box-counting dimension) was related to the Banaji--Rutar lower bounds for intermediate dimension. This connection is most naturally expressed in terms of reciprocal dimension: for sets $E \subset \R^d$, the Banaji--Rutar estimate may be conveniently formulated as
$$
\frac1{\dim_\theta E} \le \frac1\theta \, \frac1{\overline{\dim}_B E} - \left( \frac1\theta - 1 \right) \, \frac1d
$$ 
which attains an even more suggestive form when written as follows:
$$
\frac1{\dim_\theta E}  - \frac1d \le \frac1\theta \, \left( \frac1{\overline{\dim}_B E} - \frac1d \right).
$$
See \cite[Theorem 7.4]{FraTyson_intermed_dim} for an application to quasiconformal dimension distortion.
\end{remark}

\section{Applications}

In this section we discuss applications for our dimension distortion theorems. In section \ref{subsec:hardy-inequalities} we verify Hardy--Sobolev inequalities on domains obtained as quasiconformal images of domains whose boundaries enjoy good regularity properties. In section \ref{subsec:qc-classification-planar-open-sets} we distinguish cusp domains quasiconformally at the level of the dilatation. Finally, in section \ref{subsec:qc-classification-planar-open-sets-2} we provide a similar analysis for the sectorial sequence open sets defined in Example \ref{ex:another-example}.

\subsection{Hardy--Sobolev inequalities and quasiconformal mappings}\label{subsec:hardy-inequalities}

Suppose $\Omega\subset \R^d$ is open, and $u:\R^d\to \R$ is Lipschitz. We write $u\in \Lip_0(\Omega)$ if $u$ has compact support in $\Omega$.

For $x\in\Omega$, we denote the distance to the complement by
$$
d_\Omega(x) := \dist(x,\R^d\setminus\Omega).
$$ We recall the relevant weighted Hardy inequalities.

\begin{definition}\label{def:HS}
Let $\Omega\subset \R^d$ be open.
We say that $\Omega$ admits a {\it $(q,\beta)$-Hardy--Sobolev inequality} if there exists a constant
$C>0$ such that
$$
\int_{\Omega} |u(x)|^q\, d_\Omega(x)^{\beta-q}\, dx
\le
C \int_{\Omega} |\nabla u(x)|^q\, d_\Omega(x)^{\beta}\, dx,
$$
for every $u\in \Lip_0(\Omega)$.
\end{definition}

Control of the dimension of the boundary suffices to verify Hardy--Sobolev inequalities of this type. For our purposes, the formulation by Lehrb\"ack in \cite[Corollary 1.3]{HardyApplLehr} is convenient.

\begin{theorem}[Lehrb\"ack]\label{Cor: Lehrback dim Hardy}
Let $\Omega \subset \R^d$ be an open set, $1<q<\infty$, and  $\beta<q-1$.
\begin{enumerate}
\item
If $\dim_A(\R^d \setminus \Omega) < d - q + \beta$,
then $\Omega$ admits a $(q,\beta)$-Hardy inequality.

\item
If $\Omega$ is bounded, and $\dim_L(\R^d \setminus \Omega) > d - q + \beta$,
then $\Omega$ admits a $(q,\beta)$-Hardy inequality.
\item 
If $\Omega$ and $\R^d\setminus\Omega$ are both unbounded, and $\dim_L(\R^d \setminus \Omega) > d - q + \beta$,
then $\Omega$ admits a $(q,\beta)$-Hardy inequality.
\end{enumerate}
\end{theorem}

\begin{proof}[Proof of Theorem \ref{Thm: Hardy application R2}]
(a) If $\dim_A(\Omega^c) < 2-s$, then by the planar case of \cite[Theorem 1.2]{OurQCspec} we have
$$
\dim_A\bigl(f(\Omega^c)\bigr) \le \frac{p\,\dim_A(\Omega^c)}{p-2+\dim_A(\Omega^c)}
< \frac{p(2-s)}{p-2+(2-s)}.
$$
Also, by Astala's theorem \cite{Astala} we have
$$
p=\frac{2K}{K-1},
\qquad
p-2=\frac{2K}{K-1}-2=\frac{2}{K-1},
$$
which with the above imply
$$
\dim_A\bigl(f(\Omega^c)\bigr)
< \frac{\frac{2K}{K-1}(2-s)}{\frac{2}{K-1}+(2-s)}
= \frac{2K(2-s)}{2+(K-1)(2-s)}.
$$
Since
$$
\frac{2K(2-s)}{2+(K-1)(2-s)}
=2-\frac{2s}{2+(K-1)(2-s)},
$$
the conclusion follows.

\noindent (b) Noting that $f$ is a global homeomorphism of $\R^2$, impyling that $\R^2\setminus\Omega$ is unbounded if and only if $\R^2\setminus f(\Omega)$ is unbounded, the proof is similar to that of (a).
\end{proof}

\begin{remark}
Theorem \ref{Thm: Hardy application R2} cannot be obtained by change of variables, even if Definition \ref{def:HS} is relaxed to allow for $u$ to lie in some Sobolev space instead of $\operatorname{Lip}_0$. 
\end{remark}

\begin{remark}
There is a partial converse to Theorem \ref{Cor: Lehrback dim Hardy}. To wit, if a domain admits a Hardy--Sobolev inequality, then either the Hausdorff dimension of the complement is `large', or the Assouad dimension of the complement is `small'. See \cite[Theorem 6.1]{HardyApplLehr} for a precise statement. It is natural to ask whether the assumptions on Assouad and lower dimension in Theorem \ref{Thm: Hardy application R2} could be replaced with assumptions involving Hardy inequalities in the source. However, a bound on Hausdorff dimension does not yield a Hardy--Sobolev inequality, see \cite[Section 8]{HardyApplLehr}.
\end{remark}

\begin{remark}
    Working exactly as in the proof of Theorem~\ref{Thm: Hardy application R2}, we could establish similar Hardy inequality exponent estimates for domains in $\R^d$, $d\geq 3$. However, such estimates will depend on the value of the reverse H\"older inequality exponent, which is not known precisely outside the planar case. 
\end{remark}

\subsection{Dilatation-dependent quasiconformal classification of cusp domains}\label{subsec:qc-classification-planar-open-sets}

Consider the cusp domains $\Omega_\alpha$, $\alpha>1$, defined in Example \ref{ex:cusp-domains}. Ntalampekos \cite[Corollary 4.3]{nt:characterization} has proved that $\Omega_\alpha$ and $\Omega_\beta$ are quasiconformally equivalent for any $\alpha,\beta>1$. One naturally wonders what is the minimal dilatation of a quasiconformal mapping $f:\R^2 \to \R^2$ with $f(\Omega_\alpha) = \Omega_\beta$. We obtain a lower bound on the dilatation of any such map by appealing to our lower spectrum distortion theorem.

\begin{proof}[Proof of Theorem \ref{th:ntalampekos-lower-bound}]
Recall that 
\begin{equation}\label{eq:dimOmegaalpharedux}
\dim_{L,reg}^\theta(\Omega_\alpha) = \max \left\{ \frac{2-(1+\alpha)\theta}{1-\theta} , 1 \right\}, \qquad 0<\theta<1.
\end{equation}
Assume that $f:\R^2 \to \R^2$ is a homeomorphism with $f(\Omega_\alpha) = \Omega_\beta$. Set $q = \tfrac{\beta-1}{\alpha-1}$ and note that $q \in (1,\infty)$ by the assumptions on $\alpha$ and $\beta$. Denote by 
$$
K_*(q) = \left( q - \frac1{27} + \sqrt{q^2 - \frac2{27}q} \right)^{1/3} + \left( q - \frac1{27} - \sqrt{q^2 - \frac2{27}q} \right)^{1/3} - \frac13
$$
the expression on the right hand side of the inequality in \eqref{eq:Klb}. Since
$$
K_*(q) \le \sqrt{q}
$$
for all $q>0$, we may assume without loss of generality that $K(f) \le \sqrt{q}$. Then $K^2(\alpha-1) \le \beta-1$.

From \eqref{eq: dim-spec-dimension-2} we have
\begin{equation}\label{eq:planardimLspecdistortionredux}
\frac1{\dim_{L,\reg}^{\theta(t)}(\Omega_\beta)} - \frac12
\le 
K \left( \frac1{\dim_{L,\reg}^{\theta(t/K^2)}(\Omega_\alpha)} - \frac12 \right), \qquad \forall \, t > 0.
\end{equation}
We obtain a lower bound on $K$ by optimizing \eqref{eq:planardimLspecdistortionredux} over all $t>0$ using the expression in \eqref{eq:dimOmegaalpharedux}. We consider three cases according to the value of $t$.

\medskip

\noindent {\bf Case 1.} $\theta(t) > \tfrac1\beta$ and $\theta(t/K^2) > \tfrac1\alpha$, or equivalently, $t<\beta-1$ and $t/K^2<\alpha-1$. In this case, we have 
$\dim_{L,reg}^{\theta(t)}(\Omega_\beta) = \dim_{L,reg}^{\theta(t/K^2)}(\Omega_\alpha) = 1$ and \eqref{eq:planardimLspecdistortionredux} provides no information.

\medskip

\noindent {\bf Case 2.} $\theta(t) > \tfrac1\beta$ and $\theta(t/K^2) < \tfrac1\alpha$, or equivalently, $t<\beta-1$ and $t/K^2 > \alpha-1$. In this case, we have 
$\dim_{L,reg}^{\theta(t)}(\Omega_\beta) = 1$ and 
$$
\dim_{L,reg}^{\theta(t/K^2)}(\Omega_\alpha) = \frac{2-(1+\alpha)\theta(t/K^2)}{1-\theta(t/K^2)} = \frac{2t - (\alpha-1)K^2}{t}.
$$
The inequality in \eqref{eq:planardimLspecdistortionredux} reads
$$
\frac12 \le K \left( \frac{t}{2t - (\alpha-1)K^2} - \frac12 \right)
$$
which reduces to
$$
K^3 + K^2 \ge 2\frac{t}{\alpha-1}.
$$
The constraints on $t$ in this case are $K^2(\alpha-1) < t < \beta - 1$. We obtain the strongest conclusion by passing to the limit as $t \nearrow \beta-1$, namely,
\begin{equation}\label{eq:after0}
K^3 + K^2 \ge 2 \, \frac{\beta-1}{\alpha-1} = 2q.
\end{equation}

\medskip

\noindent {\bf Case 3.} $\theta(t) < \tfrac1\beta$ and $\theta(t/K^2) < \tfrac1\alpha$, or equivalently, $t > \beta-1$ and $t/K^2 > \alpha-1$.
In this case, we have 
$$
\dim_{L,reg}^{\theta(t)}(\Omega_\beta) = \frac{2-(1+\beta)\theta(t)}{1-\theta(t)} = \frac{2t-(\beta-1)}{t}
$$
and 
$$
\dim_{L,reg}^{\theta(t/K^2)}(\Omega_\alpha) = \frac{2-(1+\alpha)\theta(t/K^2)}{1-\theta(t/K^2)} = \frac{2t - (\alpha-1)K^2}{t}.
$$
The inequality in \eqref{eq:planardimLspecdistortionredux} reads
$$
\frac{t}{2t-(\beta-1)} - \frac12 \le K \left( \frac{t}{2t - (\alpha-1)K^2} - \frac12 \right)
$$
which, after some work, reduces to
\begin{equation}\label{eq:before}
K^3 \ge \frac{\beta-1}{\alpha-1} + \frac{(\beta-1)(K^3-K^2)}{2t}.
\end{equation}
We obtain the strongest conclusion by letting $t \to \beta-1$ in \eqref{eq:before}, in which case we again recover \eqref{eq:after0}.
\end{proof}

\subsection{Dilatation-dependent quasiconformal classification of special sectorial sequence open sets}\label{subsec:qc-classification-planar-open-sets-2}

Finally, we study the sectorial sequence open sets defined in Example \ref{ex:another-example}. The primary reason for defining this family of open sets is that such sets are naturally preserved by radial stretch mappings.

\begin{proposition}\label{prop:another-proposition}
Assume that $\Omega$ is a special sectorial sequence with scaling parameters $(b,a)$. Fix $K\ge 1$ and let $f:\C \to \C$ be the radial stretch map $f(z) = |z|^{1/K-1} z$. Then $\Omega' = f(\Omega)$ is a special sectorial sequence with scaling parameters $(b',a')$, where $a' = a^{1/K}$ and $b' = a^{1/K-1}b$.
\end{proposition}

\begin{proof}
Let $z \in \Omega$ and set $w = f(z)$. Then $|z| = |w|^K$. Since $a^j < |z| < a^j + \beta_j$ we have
$$
(a')^j < |w| < (a')^j + \beta_j'
$$
where
$$
\beta_j' = (a')^j \left( (1+\tfrac{\beta_j}{a^j})^{1/K}-1 \right).
$$
Moreover, $\arg(w) = \arg(z)$. We claim that there exists a constant $C_0' = C_0'(C_0,a,b,K)$ so that
$$
\tfrac1{C_0'} (b')^j \le \beta_j' \le C_0' (b')^j
$$
and
$$
\frac1{C_0'} \frac{\beta_j'}{(a')^j} \le \varphi_j \le C_0' \frac{\beta_j'}{(a')^j}.
$$
Using the elementary inequality $(1+x)^{1/K} - 1 \le \tfrac1K x$ we estimate
$$
\beta_j' \le (a')^j \frac1K \frac{\beta_j}{a^j} \le \frac{C_0}{K} \frac{(a')^j b^j}{a^j} \le C_0 (b')^j.
$$
and concavity of the function $x \mapsto (1+x)^{1/K}-1$ on the interval $[0,(\tfrac12-a)/b]$ yields
\begin{equation*}\begin{split}
\beta_j' \ge \frac{\beta_j}{a^j}
\left( \frac{ (\tfrac12-a)/b }{ (1+(\tfrac12-a)/b)^{1/K}-1 } \right)^{-1} (a')^j \gtrsim (ba'/a)^j = (b_j')
\end{split}\end{equation*}
and
$$
\varphi_j' = \varphi_j \approx (\frac{b}{a})^j = (\frac{b'}{a'})^j.
$$
This completes the proof.
\end{proof}

Now let 
\begin{equation}\label{eq:f-K-0}
f(z) = |z|^{1/K_0-1} z
\end{equation}
be the $K_0$-quasiconformal radial stretch map, and let $E$ be the closure of a special sectorial sequence $\Omega$ with scaling parameters $(b,a)$. Then $f(E)$ is the closure of a special sectorial sequence $\Omega'$ with scaling parameters $(b',a')$, $a' = a^{1/K_0}$, $b' = a^{1/K_0-1} b$. Suppose that $g:\R^2 \to \R^2$ is a $K$-quasiconformal map with $g(E) = E'$. Fix $t>0$ so that 
\begin{equation}\label{eq:theta-t}
\theta(t) = \frac{\log(1/a')}{\log(1/b')}.
\end{equation}
Then $\dim_L^{\theta(t)} \Omega' = 0$, whence also $\dim_L^{\theta(t/K^2)} \Omega = 0$, and so
\begin{equation}\label{eq:theta-t-prime}
\theta(t/K^2) \ge \frac{\log(1/a)}{\log(1/b)}.
\end{equation}
Rewriting \eqref{eq:theta-t} and \eqref{eq:theta-t-prime} we obtain
$$
t = \frac{\log(1/b')}{\log(1/a')} - 1 = 1- K_0 + K_0 \frac{\log(1/b)}{\log(1/a)} - 1 = K_0 \left(  \frac{\log(1/b)}{\log(1/a)} - 1 \right)
$$
and
$$
\frac{t}{K^2} \le \frac{\log(1/b)}{\log(1/a)} - 1.
$$
Hence $K_0/K^2 \le 1$ and so $K \ge \sqrt{K_0}$. In conclusion, we have shown the following result.

\begin{proposition}\label{prop:ssos-lower-bounds}
Let $\Omega$ and $\Omega'$ are special sectorial sequences with scaling parameters $(b,a)$ and $(b',a')$ respectively. Then the radial stretch map $f$ as in \eqref{eq:f-K-0} with
$$
K_0 = \frac{\tfrac{\log(1/b')}{\log(1/a')} - 1}{\tfrac{\log(1/b)}{\log(1/a)} - 1}
$$
is a $K_0$-quasiconformal map with $f(\Omega) = \Omega'$. Moreover, if $g:\R^2 \to \R^2$ is any $K$-quasiconformal map with $g(\Omega) = \Omega'$, then $K \ge \sqrt{K_0}$.
\end{proposition}

\begin{remark}
We conclude this section by stating a general proposition which encapsulates the principle at work in Theorem \ref{th:ntalampekos-lower-bound} and Proposition \ref{prop:ssos-lower-bounds}. The proof is elementary and is left to the reader. Note that we assume that $\Omega$ and $\Omega'$ are homeomorphic open sets, but we do not assume that they are domains (i.e., connected).
\end{remark}

\begin{proposition}\label{prop:qc-classification-by-lower-spectrum}
Let $\Omega,\Omega'$ be homeomorphic open sets of $\R^d$ and let $K_\ddag(\Omega,\Omega')$ be the infimum of values $K \ge 1$ for which there exists a $K$-quasiconformal map $f:\R^d \to \R^d$ with $f(\Omega) = \Omega'$. 

For $0<s<d$ set $\Phi(s) = s^{-1}-d^{-1}$, and for $p>d$ set $\alpha(p) = 1-d/p$. For each $t>0$, let $K(t)$ denote the infimum of values $K\ge 1$ so that
$$
\Phi(\dim_{L,\reg}^{\theta(t)}(\Omega')) \le
\alpha(p_I^{\RH}(d,K)) \, \Phi(\dim_{L,\reg}^{\theta(\tilde{t}/K^{d-1})}(\Omega')).
$$
Then
$$
K_\ddag(\Omega,\Omega') \} \ge \sup_{t>0} K(t).
$$
\end{proposition}

\section{Quasiconformal and Sobolev distortion of the lower box dimension}\label{sec:lowerboxdim}

Recall that for a bounded set $E\subset \R^d$, the lower box dimension of $E$ is defined by
$$
\underline{\dim}_B E
=
\liminf_{r\rightarrow0}
\frac{\log N_r(E)}{-\log r}.
$$
Following similar methods to those in the proof of Theorem~\ref{thm: QC L spectrum distortion} yields lower box dimension distortion bounds for local Sobolev mappings. This naturally yields similar bounds for quasiconformal mappings, with the improved Sobolev higher integrability exponent instead of the reverse H\"older inequality exponent. Although the set-up is slightly different, this result is the closest to and complements that of Kaufman in \cite{Kaufman} on the upper box dimension.

We first characterize the lower box dimension in terms of a growth condition on ball covering numbers, similar to that in the definition of lower dimension.

\begin{proposition}\label{prop: lower-box-dim-spectrum}
Let $E\subset\R^d$ be a non-empty bounded set. Then
$$
\underline{\dim}_B E
=
\sup\left\{
s\geq0:
\exists\, C>0\text{ such that }
N_r(E)\geq (1/C) r^{-s}
\text{ for all }0<r<\diam E
\right\}.
$$
\end{proposition}

\begin{proof}
Denote the supremum on the right-hand side of the relation by $D$.

Suppose $s\geq0$ is such that there exists $C>0$ satisfying
$$
N_r(E)\geq (1/C) r^{-s}
$$
for all $0<r<\diam E$. Taking logarithms gives
$$
\frac{\log N_r(E)}{-\log r}
\geq
s-\frac{\log C}{-\log r}.
$$
Taking $\liminf_{r\to 0}$ in the above yields
$$
\underline{\dim}_B E
=
\liminf_{r\rightarrow0}
\frac{\log N_r(E)}{-\log r}
\geq s.
$$
Since $s$ is arbitrary, taking the supremum gives $D\leq\underline{\dim}_B E$.

For the other direction, let
$$
0\leq s<\underline{\dim}_B E,
$$
and choose $u$ such that
$$
s<u<\underline{\dim}_B E.
$$
By the definition of the limit inferior, there exists $r_0\in(0,1)$ such that
$$
\frac{\log N_r(E)}{-\log r}\geq u
$$
for every $0<r<r_0$. Consequently, $N_r(E)\geq r^{-u}\geq r^{-s}$ for every $0<r<r_0$.

On the other hand, for $r_0\leq r<\diam E$ we have
$$
N_r(E)\geq1
\geq
r_0^s r^{-s},
$$
which implies $s\leq D$. Letting
$s\rightarrow\underline{\dim}_B E$ gives
$$
\underline{\dim}_B E\leq D,
$$
which completes the proof.
\end{proof}

This characterization of the lower box dimension allows us to follow a similar route to the proof of Theorem~\ref{thm: QC L spectrum distortion}. We begin with an upper estimate valid for arbitrary supercritical Sobolev mappings.

\begin{theorem}\label{thm: Sob ldimB}
Let $f\in W^{1,p}_\loc(\Omega:\R^N)$ continuous, for a domain $\Omega \subset \R^d$, with $p>d$ and $N \in \N$. If $E\subset\Omega$ is a compact set with $\alpha_E=\underline{\dim}_B E$, then
$$
\underline{\dim}_B f(E)
\leq
\frac{p\,\alpha_E}
{p-d+\alpha_E}.
$$
\end{theorem}

\begin{proof}
After composing $f$ with appropriate similarities, we may assume that 
$$
\max\{\diam E, \diam f(E)\}\leq 1/10,
$$ 
and proceed as in the proof of Theorem \ref{thm: QC L spectrum distortion}, making use the equivalent definition of the lower box dimension in Proposition \ref{prop: lower-box-dim-spectrum}. The arguments are simpler in this case, due to the large scales $R_n$ and $R_n'$ being replaced by $\diam E$ and $\diam f(E)$, respectively, and the two spectrum variables $\theta(Kt)$ and $\theta((1-d/p)t)$ being replaced by the constant $1$. The arguments involving major cubes and the Morrey--Sobolev inequality remain  unchanged, with the use of the reverse H\"older inequality being unnecessary in this situation;
the convergence of
$$
\int_{\bigcup\limits_{k\geq m}\,\bigcup\limits_{i} Q_i^k} |Df|^p
$$ 
follows directly by absolute continuity of the integral. 
\end{proof}

Recalling the higher integrability of $K$-quasiconformal mappings, the following corollary naturally follows.

\begin{corollary}\label{cor: QC ldimB}
    Let $f:\Omega\rightarrow\Omega'$ be a $K$-quasiconformal mapping between domains in $\R^d$, and let $E\subset\Omega$ be a non-empty compact set with $\alpha_E=\underline{\dim}_B E\in (0,d)$.
    Then 
    $$
    \frac{(p_I^{\Sob}(d,K)-d)\, \alpha_E}{p_I^{\Sob}(d,K)-\alpha_E}
    \le
    \underline{\dim}_B f(E)
    \leq
    \frac{p_O^{\Sob}(d,K)\,\alpha_E}
    {p_O^{\Sob}(d,K)-d+\alpha_E}.
    $$
\end{corollary}

\begin{proof}
    The right-hand side follows by Theorem~\ref{thm: Sob ldimB}, and the left-hand side follows by applying the right-hand side inequality to $f^{-1}$, along with $K_O(f^{-1})=K_I(f)$.
\end{proof}

\bibliographystyle{acm}
\bibliography{QCLowerDim}
\end{document}